\documentclass[12pt]{article}  

\usepackage{amsfonts}
\usepackage{amssymb}
\usepackage[all]{xy}

\newcommand{\cA}{{\cal A}}
\newcommand{\cC}{{\cal C}}
\newcommand{\cB}{{\cal B}}

\newcommand{\cF}{{\cal F}}
\newcommand{\cG}{{\cal G}}
\newcommand{\cH}{{\cal H}}

\newcommand{\cK}{{\cal K}}
\newcommand{\cL}{{\cal L}}
\newcommand{\cM}{{\cal M}}

\newcommand{\cO}{{\cal O}}
\newcommand{\cP}{{\cal P}}
\newcommand{\cQ}{{\cal Q}}

\newcommand{\cV}{{\cal V}}

\newcommand{\gS}{{\Sigma}}
\newcommand{\gL}{{\Lambda}}
\newcommand{\gD}{{\Delta}}

\newcommand{\gs}{{\sigma}}
\newcommand{\gd}{{\delta}}
\newcommand{\gG}{{\mit\Gamma}}

\newcommand{\Sym}{{\rm Sym}}

\newcommand{\op}{{\rm op}}
\newcommand{\coker}{{\rm coker}}

\newcommand{\End}{{\rm End}}
\newcommand{\Hom}{{\rm Hom}}

\newcommand{\id}{{\rm id}}
\newcommand{\Spec}{{\rm Spec}}

\newcommand{\Ind}{{\rm Ind}}
\newcommand{\Res}{{\rm Res}}
\newcommand{\sgn}{{\rm sgn}}
\newcommand{\Tor}{{\rm Tor}}
\newcommand{\rmcr}{{\rm cr}}
\newcommand{\Image}{{\rm Im}}
\newcommand{\Kos}{{\rm Kos}}
\newcommand{\Simp}{{\rm Simp}}
\newcommand{\Compl}{{\rm Compl}}
\newcommand{\Sur}{{\rm Sur}}
\newcommand{\Tot}{{\rm Tot}}

\newcommand{\ZZ}{{\mathbb Z}}

\newcommand{\NN}{{\mathbb N}}

\newcommand{\QQ}{{\mathbb Q}}

\newcommand{\ra}{\rightarrow}
\def\rightepi{{\longrightarrow \kern-0.7em \rightarrow}}
\newcommand{\oplusm}{\mathop{\oplus}\limits}

\newcommand{\notteilt}{{\,\not{\kern-0.075em|}\,}}
\def\antiddots{\mathinner{\mkern1mu\raise1pt\vbox{\kern7pt\hbox{.}}\mkern2mu
    \raise4pt\hbox{.}\mkern2mu\raise7pt\hbox{.}\mkern1mu}}

\begin{document}

\vspace*{15ex}

\begin{center}
{\LARGE\bf Computing the Homology of Koszul Complexes}\\
\bigskip
by\\
\bigskip
{\sc Bernhard K\"ock}
\end{center}

\bigskip

\begin{quote}
{\footnotesize {\bf Abstract}. Let $R$ be a commutative ring and $I$ an ideal in $R$ which
is locally generated by a regular sequence of length $d$. Then, each f.\ g.\ projective 
$R/I$-module $V$ has an $R$-projective resolution $P.$ of length $d$. In this paper, we 
compute the homology of the $n$-th Koszul complex associated with the homomorphism 
$P_1 \ra P_0$ for all $n \ge 1$, if $d=1$. This computation yields a new proof of the
classical Adams-Riemann-Roch formula for regular closed immersions which does not use
the deformation to the normal cone any longer. Furthermore, if $d=2$, we compute the 
homology of the complex $N\, \Sym^2 \, \gG(P.)$ where $\gG$ and $N$ denote the functors 
occurring in the Dold-Kan correspondence.}
\end{quote}

\bigskip

\section*{Introduction}

In the theory of exterior powers of f.\ g.\ projective modules over a commutative
ring $R$, there are various formulas which, so far, can only be proved by using
geometric means (e.g., the projective fibre bundle or blowing up). An example
is the formula which computes the composition of exterior powers in the 
Grothendieck ring $K_0(R)$ (e.g., see \cite{FL}). The plethysm problem is to
find (explicit) functorial short exact sequences which
prove this formula within the framework of commutative algebra. 
It has been intensively studied in universal
representation theory. In this paper, we study the analogous problem for the
Adams-Riemann-Roch formula for regular closed immersions (between affine schemes).

We now recall the Adams-Riemann-Roch formula (see \cite{FL}). Let $i: 
Y \hookrightarrow X$ be a regular closed immersion of schemes with the locally free
conormal sheaf $\cC$. Let $\psi_n : K_0(X) \ra K_0(X)$ denote the $n$-th
Adams operation on the Grothendieck ring $K_0(X)$; it is defined as a certain
integral polynomial in the exterior power operations on $K_0(X)$. For any locally
free $\cO_Y$-module $\cV$, we choose a locally free resolution $\cP.$
of $i_*(\cV)$ on $X$, and we map the element $[\cV] \in K_0(Y)$ to the Euler
characteristic $\sum_{k \ge 0} (-1)^k [\cP_k] \in K_0(X)$. This association
induces a well-defined homomorphism $i_*: K_0(Y) \ra K_0(X)$. The 
Adams-Riemann-Roch theorem for $i_*$ now states that, for all $y \in K_0(Y)$, we have:
\[\psi_n(i_*(y)) = i_*(\theta^n(\cC) \cdot \psi_n(y)) \quad {\rm in} \quad K_0(X);\]
here, $\theta^n(\cC)$ denotes the so-called Bott element associated with $\cC$ 
(see \cite{FL}). 

This formula (for one $n\ge 2$) implies the famous Grothendieck-Riemann-Roch theorem
for $i_*$ in a rather formal manner (see \cite{FL}). It has been formulated for the
first time by Manin in \cite{Ma}, but it has been proved there only up to 
$n^2$-torsion in $K_0(X)$. Moreover, as a formula in $K_0(X) \otimes \QQ$, 
it follows already
from Th\'eor\`eme 4.3 in Expose VII of \cite{SGA6}. The given version without
denominators can be derived from Th\'eor\`eme 2.1 in Jouanolou's paper \cite{J}
and can be found in Soul\'e's paper \cite{So} and in the book \cite{FL} by
Fulton and Lang. In all these sources, the deformation to the normal cone is
the crucial tool in the proof. The object of 
this paper is to construct explicit short exact sequences of $\cO_X$-modules 
which prove this formula
in the case ${\rm codim}(Y/X) = 1$, $n$ arbitrary and in the case ${\rm codim}(Y/X)
= 2$, $n=2$. Since the proof of the Adams-Riemann-Roch formula in the general case
can immediately be reduced to the case ${\rm codim}(Y/X) =1$ using blowing up and
the excess intersection formula (see \cite{SGA6}), we in particular obtain a new,
very natural and simple proof in the general case which, in contrast to \cite{J},
\cite{So}, or \cite{FL}, does not use the deformation any longer.

To illustrate our approach, we now
rather completely describe it for ${\rm codim}(Y/X) =1$ and
$n=2$. As in the beginning of the introduction, we consider only the affine case;
so, let $X=\Spec(R)$ and $Y = \Spec(R/I)$ where the ideal $I$ of $R$ is locally
generated by a non-zero-divisor. Let $\sigma_2$ and $\lambda_2$ denote the
second symmetric power and second exterior power operation. Let $V$ be a f.\ g.\
projective $R/I$-module and 
\[0 \ra P \ra Q \ra V \ra 0\]
an $R$-projective resolution of $V$. Since $\psi_2 = \sigma_2 - \lambda_2$ and
$\theta^2(I/I^2) = 1 + [I/I^2]$, it suffices to prove the formula
\begin{equation}
\sigma_2(i_*([V])) = i_*([\Sym^2(V)] - [\gL^2(V) \otimes I/I^2]) \quad 
{\rm in} \quad K_0(X)
\end{equation}
and a similar formula for $\lambda_2$. Now, the main idea is to consider the Koszul
complex
\[ K: \quad 0 \ra \gL^2(P) \ra P\otimes Q \ra \Sym^2(Q) \ra 0\]
which, for instance, can be defined as a subcomplex of the natural short exact 
sequence $0 \ra \gL^2(Q) \ra Q\otimes Q \ra \Sym^2(Q) \ra 0$. It is 
immediately clear that the alternating sum of the objects of this complex 
equals the left hand side of the formula (1). Since the alternating sum
of the objects equals the alternating sum of the homology modules, it suffices
to show that $H_0(K) \cong \Sym^2(V)$, $H_1(K) \cong \gL^2(V) \otimes I/I^2$, and
$H_2(K) \cong 0$. (Here, $H_0(K)$, $H_1(K)$, and $H_2(K)$ denote the homology of
$K$ at the places $\Sym^2(Q)$, $P\otimes Q$, and $\gL^2(P)$, respectively.) Whilst
the isomorphisms $H_0(K) \cong \Sym^2(V)$ and $H_2(K) \cong 0$ can rather easily
be proved (see Remark 3.3), the proof of the isomorphism $H_1(K) \cong
\gL^2(V) \otimes I/I^2$ is more complicated. The most important observation for this
is the fact that $H_1(K)$ does not depend on the chosen resolution and that
the association $V \mapsto H_1(K) =: F(V)$ is even a functor. This can be proved
as follows. Any homomorphism $V \ra V'$ between f.\ g.\ projective $R/I$-modules
can be lifted to a homomorphism from a resolution of $V$ to a resolution of $V'$,
and this lifted homomorphism then induces a homomorphism between the corresponding
Koszul complexes $K$ and $K'$ in the obvious way. In order 
to prove that the induced homomorphism
between the homology modules does not depend on the lifting, we consider
the complex 
\[ L: \quad 0 \ra P\otimes P \ra \Sym^2(P) \oplus P\otimes Q \ra 
\Sym^2(Q) \ra 0\]
(with the obvious differentials) which is quasi-isomorphic to the complex $K$. 
One can show (see also Remark 2.6) that any
homotopy between two homomorphisms between the resolutions 
induces a homotopy between
the corresponding homomorphisms between the complexes $L$ and $L'$. (The analogous
statement for the complexes $K$ and $K'$ is not true!) As usual in homological 
algebra, this implies that the association $V \mapsto H_1(L) \cong F(V)$ is
a well-defined functor. In particular, for any f.\ g.\ projective $R/I$-modules
$V,W$, the direct sum $F(V) \oplus F(W)$ is a canonical direct summand of $F(V \oplus W)$.
Let $\rmcr_2(F)(V,W)$ denote the complement. Then $\rmcr_2(F)$ is a bifunctor. It is 
called the second cross effect of $F$. (This terminology has been introduced by
Eilenberg and MacLane in \cite{EM}.) One easily sees that we have $\rmcr_2(F) \cong
\Tor_1^R$ in our situation. Hence, $\rmcr_2(F)$ is biadditive, and we have
$\rmcr_2(F)(R/I,R/I) \cong I/I^2$. Furthermore, one immediately sees that 
$F(R/I) \cong 0$. Eilenberg (\cite{E}) and Watts (\cite{Wa}) have proved that
any additive functor $F$ on the category of f.\ g.\ projective modules over a 
commutative ring $A$ is isomorphic to the tensor functor $F(A) \otimes -$. 
Similarly, one easily shows that any functor $F$ with the three properties just
proved is isomorphic to the functor $\gL^2(-) \otimes I/I^2$ (see Corollary 1.6). 
This finally proves the claimed isomorphism $H_1(K) \cong \gL^2(V) \otimes I/I^2$
and hence the formula (1). The analogous formula for $\lambda_2$ can be proved
in a similar manner (see also Remark 3.6).

Now we briefly describe the idea of our approach in the case ${\rm codim}(Y/X) =2$,
$n=2$. Let $\cV$ be a locally free $\cO_Y$-module and $\cP.$ a locally free
resolution of $i_*(\cV)$ on $X$. Here, we compute the homology of the complex
$\cL := N \Sym^2 \gG(\cP.)$ where $\gG$ and $N$ denote the 
functors occurring in the Dold-Kan correspondence (see section 2). The idea to 
consider this complex is suggested by the observation that the complex $L$
used in the case ${\rm codim}(Y/X) =1$ is isomorphic to the complex 
$N \Sym^2 \gG(P\ra Q)$ (see section 2). On the contrary, it seems that an analogue of the 
Koszul complex $K$ does not exist in the case ${\rm codim}(Y/X) \ge 2$. We
obtain the following result (see Theorem 6.4):
$H_0(\cL) \cong i_*(\Sym^2(\cV))$, $H_1(\cL) \cong i_*(\gL^2(\cV) \otimes \cC)$,
$H_2(\cL) \cong i_*(D^2(\cV) \otimes \gL^2(\cC))$, and $H_k(\cL) \cong 0 \,\, 
{\rm for} \,\, k \ge 3$.
This computation of the homology then implies the Adams-Riemann-Roch formula
in the case ${\rm codim}(Y/X) = 2 $, $n=2$ as above (see section 6). The 
main tools we use in this computation are the Eilenberg-Zilber theorem
and the universal form of the Cauchy decomposition developed by Akin,
Buchsbaum and Weyman in \cite{ABW}. We again verify certain properties of
the corresponding cross effect functors which, together with easy abstract
lemmas, imply the above isomorphisms.

In section 5, we present a further application of the theory of cross effect functors
developed in section 1. We prove the Riemann-Roch formula for tensor power operations
(see Theorem 5.3) which has been proved in the paper \cite{KoTe} using the 
deformation to the normal cone (and which in turn implies the Adams-Riemann-Roch
formula in characteristic $0$, see section 4 in \cite{KoTe}). 
Here, the codimension of $Y$ in $X$ may be arbitrary.
For any locally free $\cO_Y$-module $\cV$ and for any locally free resolution 
$\cP. \ra i_*(\cV)$ of $i_*(\cV)$ on $X$, we compute the homology of the total
complex of the $n$-th tensor power $\cP.^{\otimes n}$ of $\cP.$ 
together with the canonical action of the symmetric group (see Theorem 5.1).

{\bf Acknowledgments.} I would like to thank D.\ Grayson and R.\ McCarthy
for many very helpful discussions where I have learned, for example, the
(advantages of) simplicial techniques. In particular, D.\ Grayson has 
contributed the basic idea for the important Corollary 2.5. Furthermore, R.\ McCarthy has
explained the main features of cross effect functors to me. It is a pleasure
to thank them for all this assistance.

\bigskip

\section*{\S 1 On Cross Effect Functors}

First, we recall the definition of the cross effect functors $\rmcr_k(F)$, $k\ge 0$,
associated with a functor $F$ between additive categories (cf.\ \cite{EM} and \cite{JM}).
Then, we introduce the diagonal and plus maps and describe these maps explicitly for 
the second and third symmetric power functor. The most important result of this 
section is Theorem 1.5 which states that any functor $F$ (of finite degree) from the
category of f.\ g.\ projective modules over a ring $A$ to any abelian category is
essentially determined by the objects $\rmcr_k(F)(A, \ldots, A)$, $k\ge 0$, and
the diagonal and plus maps between these objects. Finally, we characterize the
classical exterior power, symmetric power and divided power functors.

Let $\cP$ be an additive category, $\cM$ an abelian category, and $F: \cP \ra \cM$ 
a functor with $F(0) = 0$. 

{\bf Definition 1.1}. Let $k \ge 0$. For any $V_1, \ldots, V_k \in \cP$ and $1 \le i \le k$,
let 
\[p_i: V_1 \oplus \ldots \oplus V_k \ra V_i \ra V_1 \oplus \ldots \oplus V_k\] 
denote the $i$-th projection. The $k$-functor
\[\begin{array}{cccl}
\rmcr_k(F): & \cP^k & \ra & \cM \\
& (V_1, \ldots, V_k) & \mapsto & \Image\left( \sum_{j=1}^k (-1)^{k-j}
\sum_{1 \le i_1 < \ldots < i_j \le k} F( p_{i_1} + \ldots + p_{i_j})\right)
\end{array}\]
is called the {\em $k$-th cross effect of $F$}; 
here, for any $V_1, \ldots,  V_k, W_1,\ldots,
W_k \in \cP$ and $f_1 \in \Hom_\cP(V_1,W_1), \:\ldots,\: f_k \in \Hom_\cP(V_k, W_k)$,
the map
\[\rmcr_k(F)(f_1, \ldots, f_k): \rmcr_k(F)(V_1, \ldots, V_k) \ra
\rmcr_k(F)(W_1, \ldots, W_k)\]
is induced by $f_1\oplus \ldots \oplus f_k \in \Hom_\cP(V_1 \oplus \ldots \oplus
V_k, W_1 \oplus \ldots \oplus W_k)$. $F$ is said to be a {\em functor of degree $\le k$}
if $\rmcr_{k+1}(F)$ is identically zero.

The cross effect functors $\rmcr_k(F)$, $k\ge 0$, have the following properties. 
We obviously have $\rmcr_0(F) \equiv 0$ and $\rmcr_1(F) = F$. Furthermore, we have
$\rmcr_k(F)(V_1, \ldots, V_k) =0$, if $V_i =0$ for one $i\in \{1, \ldots, k\}$ (see
Theorem 9.2 on p.\ 79 in \cite{EM}). The canonical action of the symmetric group
$\gS_k$ on $V^k = V \oplus \ldots \oplus V$ induces a natural action of $\gS_k$ on
$\rmcr_k(F)(V, \ldots, V)$ (see Theorem 9.3 on p.\ 80 in \cite{EM}). Any natural
transformation $F \ra G$ between functors $F$ and $G$ from $\cP$ to $\cM$ induces
a natural transformation $\rmcr_k(F) \ra \rmcr_k(G)$ (see Theorem 9.5 on p.\ 80 in
\cite{EM}). The most important property is the following proposition.

{\bf Proposition 1.2.} For any $k,l \ge 1$ and $V_1, \ldots, V_l \in  \cP$, we
have a canonical isomorphism
\[\rmcr_k(F)(\ldots,\: V_1\oplus \ldots \oplus V_l,\: \ldots) \cong
\oplusm_{1\le j \le l} \oplusm_{1\le i_1 < \ldots < i_j \le l}
\rmcr_{k+j}(F)(\ldots,\: V_{i_1}, \ldots, V_{i_j},\: \ldots)\]
which is functorial in $V_1, \ldots, V_l$. In particular, $F$ is of degree $\le k$, 
if and only if $\rmcr_k(F)$ is a $k$-additive functor. 

{\bf Proof.} See Theorem 9.1, Lemma 9.8 and Lemma 9.9 in \cite{EM}.

For instance, Proposition 1.2  states that 
\[F(V_1 \oplus V_2) \cong F(V_1) \oplus F(V_2) \oplus \rmcr_2(F)(V_1, V_2)\]
for all $V_1, V_2 \in \cP$; i.e., $\rmcr_2(F)$ measures the deviation from additivity
of the functor $F$. This isomorphism can also be used to define $\rmcr_2(F)$ (see
section 3 in \cite{JM}). Similarly, the higher cross effects can inductively be
defined by the isomorphism
\begin{eqnarray*}
\lefteqn{\rmcr_k(F)(V_1, \ldots, V_{k-1}, V_k \oplus V_{k+1}) \cong}\\
&&\rmcr_k(F)(V_1, \ldots, V_{k-1}, V_k) \oplus \rmcr_k(F)(V_1, \ldots, V_{k-1}, V_{k+1})
\oplus \rmcr_{k+1}(F)(V_1, \ldots, V_k, V_{k+1})
\end{eqnarray*}
(see section 7 in \cite{JM}). In fact, we will use this definition whenever we have
to compute cross effect functors. 

{\bf Definition 1.3.} Let $l \ge k \ge 1$, $\varepsilon = (\varepsilon_1, \ldots,
\varepsilon_k) \in \{1, \ldots, l\}^k$ with $|\varepsilon| = \varepsilon_1 
+ \ldots + \varepsilon_k = l$ and $V_1, \ldots, V_k \in \cP$. The composition
\begin{eqnarray*}
\lefteqn{\xymatrix{
\gD_\varepsilon: \rmcr_k(F)(V_1, \ldots, V_k) \ar[rrr]^-{\rmcr_k(F)(\gD, \ldots, \gD)}
&&& \rmcr_k(F)(V_1^{\varepsilon_1}, \ldots, V_k^{\varepsilon_k})} } \hspace{5cm}\\ 
&&\xymatrix{
{}\ar@{>>}[r]^-{\pi}& \rmcr_l(F)(V_1, \ldots, V_1, \:\ldots,\: V_k, \ldots, V_k) }
\end{eqnarray*}
of the map $\rmcr_k(F)(\gD, \ldots, \gD)$ (induced by the diagonal maps
$\gD: V_i \ra V_i^{\varepsilon_i}$, $i=1, \ldots, k$) with the canonical
projection $\pi$ (according to Proposition 1.2) is called {\em diagonal
map associated with $\varepsilon$}. The analogous composition
\begin{eqnarray*}
\lefteqn{ \xymatrix{
+_\varepsilon:  \rmcr_l(F)(V_1, \ldots, V_1,\: \ldots,\: V_k, \ldots, V_k)
\ar@{^{(}->}[r] & 
\rmcr_k(F)(V_1^{\varepsilon_1}, \ldots, V_k^{\varepsilon_k}) }} \hspace{7cm}\\
&& \xymatrix{\ar[rrr]^-{\rmcr_k(F)(+, \ldots, +)} &&& \rmcr_k(F)(V_1, \ldots, V_k) }
\end{eqnarray*}
is called {\em plus map associated with $\varepsilon$}.

Obviously, the maps $\gD_\varepsilon$ and $+_\varepsilon$ form natural transformations
between the functors $\rmcr_k(F)$ and $\rmcr_l(F) \circ (\gD_{\varepsilon_1}, \ldots,
\gD_{\varepsilon_k})$ from $\cP^k$ to $\cM$. One easily sees that the map
$\gD_{\varepsilon}$ can be decomposed into a composition of maps $\gD_\delta$
with $\delta \in \{1,2\}^j$ such that $|\delta| = j+1$ and $j \in \{k, \ldots, l-1\}$.
The same holds for $+_\varepsilon$.

{\bf Example 1.4.} Let $\cP = \cM$ be the category of modules over a commutative 
ring $A$. For any $n \ge 1$, let $\Sym^n: \cP \ra \cM$ denote the $n$-th 
symmetric power functor.\\
(a) For all $k \ge 1$ and $V_1, \ldots, V_k \in \cP$, we have:
\[\rmcr_k(\Sym^n)(V_1, \ldots, V_k) \cong 
\oplusm_{{(n_1, \ldots, n_k) \in \{1, \ldots, n\}^k \atop n_1+\ldots +n_k =n}}
\Sym^{n_1}(V_1) \otimes \ldots \otimes \Sym^{n_k}(V_k).\]
(b) Let $n=2$. Then the following diagrams commute for all $V\in \cP$:
\[\xymatrix@R=0pt{
&\rmcr_2(\Sym^2)(V,V) \\
\\
&{}\ar@{=}[dd]\\
&{}\\
&{}\\
\Sym^2(V) \ar[uuuuur]^-{\gD_2} \ar[r] & V\otimes V \\
v_1 v_2  \ar@{|->}[r] & v_1 \otimes v_2 + v_2 \otimes v_1
}
\quad \xymatrix@R=0pt{ \\ \\ \\ {\rm and} } \quad
\xymatrix@R=0pt{
\rmcr_2(\Sym^2)(V,V) \ar[dddddr]^-{+_2}\\
\\
{}\ar@{=}[dd]\\
\\
\\
V\otimes V \ar[r] & \Sym^2(V) \\
v_1 \otimes v_2 \ar@{|->}[r] & v_1 v_2. 
}\]
(c) Let $n =3$. By virtue of the isomorphisms given in (a), the diagonal and plus
maps can be described as follows (for $V,W \in \cP$):
\[ \begin{array}{c}
\gD_3: \Sym^3(V) \ra V \otimes V \otimes V, \quad v_1v_2v_3 \mapsto
\sum_{\sigma \in \gS_3} v_{\sigma(1)} \otimes v_{\sigma(2)} \otimes v_{\sigma(3)}.
\end{array}\]
\[\begin{array}{cccc}
\gD_{(2,1)}: &\Sym^2(V) \otimes W \oplus V\otimes \Sym^2(W) &\ra &
V\otimes V \otimes W\\
&(v_1v_2 \otimes w, v \otimes w_1w_2)& \mapsto & v_1 \otimes v_2 \otimes w +
v_2 \otimes v_1 \otimes w.
\end{array}\]
\[\begin{array}{cccc}
\gD_{(1,2)}: &\Sym^2(V) \otimes W \oplus V\otimes \Sym^2(W) &\ra &
V\otimes W \otimes W \\
&(v_1v_2 \otimes w, v\otimes w_1w_2) & \mapsto & v\otimes w_1 \otimes w_2 +
v\otimes w_2 \otimes w_1.
\end{array}\]
\[\begin{array}{c}
\gD_2:  \Sym^3(V)  \ra  \Sym^2(V) \otimes V \oplus V\otimes \Sym^2(V)\\
v_1v_2v_3  \mapsto  (v_2v_3 \otimes v_1 + v_1v_3\otimes v_2 + 
v_1v_2 \otimes v_3, v_1 \otimes v_2v_3 + v_2 \otimes v_1v_3 + v_3 \otimes v_1v_2).
\end{array}\]
\[+_3 : V\otimes V \otimes V \ra \Sym^3(V), \quad v_1\otimes v_2 \otimes v_3 
\mapsto v_1v_2v_3.\]
\[+_{(2,1)}: V\otimes V \otimes W \ra \Sym^2(V) \otimes W \oplus V \otimes \Sym^2(W),
\quad v_1 \otimes v_2 \otimes w \mapsto (v_1v_2 \otimes w, 0).\]
\[+_{(1,2)}: V \otimes W \otimes W \ra \Sym^2(V) \otimes W \oplus V \otimes \Sym^2(W),
\quad v \otimes w_1 \otimes w_2 \mapsto (0, v\otimes w_1 w_2).\]
\[+_2:  \Sym^2(V) \otimes V \oplus V\otimes \Sym^2(V) \ra \Sym^3(V), 
\quad (v_1 v_2\otimes v_3, v'_1\otimes v'_2 v'_3) \mapsto v_1v_2v_3 + v'_1v'_2v'_3.\]

{\bf Proof.} Straightforward.

The following theorem states that any functor $F$ of finite degree from the category
of all f.\ g.\ projective modules over a ring $A$ to any abelian category $\cM$ 
is determined by
the objects $\rmcr_i(F)(A, \ldots, A)$, $i\ge 0$, the diagonal and plus maps between
these objects, and the actions
\[A \ra \End_\cM(\rmcr_i(F)(A, \ldots, A)) , \quad a \mapsto 
\rmcr_i(F)(1, \ldots, 1, a, 1, \ldots, 1), \]
of the multiplicative monoid $A$ on these objects.

{\bf Theorem 1.5.} Let $A$ be a ring, $\cM$ an abelian category, $d\in \NN$, and
\[\xymatrix{
F,G: (\mbox{f.\ g.\ projective }A\mbox{-modules}) 
\ar@<0.5ex>[r] \ar@<-0.5ex>[r] & \cM }\]
two functors of degree $\le d$ with $F(0) =0= G(0)$. Suppose that there exist
isomorphisms
\[\xymatrix@1{
\alpha_i(A, \ldots, A):\rmcr_i(F)(A, \ldots, A) \ar[r]^-{\sim}&
\rmcr_i(G)(A, \ldots, A),} \quad i=1, \ldots, d,\]
which are compatible with the action of $A$ in each component and which make the following
diagrams commute for all $i \in \{1, \ldots, d-1\}$ and
$\varepsilon \in \{1,2\}^i$ with $|\varepsilon|=i+1$:
\[\xymatrix{
\rmcr_i(F)(A, \ldots, A) \ar[r]^-{\sim} \ar[d]^-{\gD_\varepsilon}
& \rmcr_i(G)(A, \ldots, A) \ar[d]^-{\gD_\varepsilon}\\
\rmcr_{i+1}(F)(A, \ldots, A) \ar[r]^-{\sim} & \rmcr_{i+1}(G)(A, \ldots, A)
}\]
and 
\[\xymatrix{
\rmcr_{i+1}(F)(A, \ldots, A) \ar[r]^-{\sim} \ar[d]^-{+_\varepsilon}& 
\rmcr_{i+1}(G)(A, \ldots, A) \ar[d]^-{+_\varepsilon}\\
\rmcr_{i}(F)(A, \ldots, A) \ar[r]^-{\sim} & \rmcr_{i}(G)(A, \ldots, A).
}\]
Then the functors $F$ and $G$ are isomorphic.

{\bf Proof.} For any $i \ge 1$ and $n_1, \ldots, n_i \ge 1$, we write 
$\rmcr_i(F)(A^{n_1}, \ldots, A^{n_i})$ as a direct sum of copies of
$\rmcr_i(F)(A, \ldots, A),\: \ldots,\: \rmcr_d(F)(A, \ldots, A)$ according to
Proposition 1.2; in the same way, we write $\rmcr_i(G)(A^{n_1}, \ldots, A^{n_i})$
as a direct sum of copies of $\rmcr_i(G)(A, \ldots, A), \:\ldots, \:
\rmcr_d(G)(A, \ldots, A)$. The corresponding sum of copies of the isomorphisms
$\alpha_i(A, \ldots, A),\: \ldots,\: \alpha_d(A, \ldots, A)$ then yields an
isomorphism
\[\xymatrix{
\alpha_i(A^{n_1}, \ldots, A^{n_i}): \rmcr_i(F)(A^{n_1}, \ldots, A^{n_i}) 
\ar[r]^-{\sim} & \rmcr_i(G)(A^{n_1}, \ldots, A^{n_i}).}\]
Now, we show by descending induction on $i$ that, for all $m_1, \ldots, m_i ,
n_1, \ldots, n_i \ge 1$ and for all matrices $\cA_1 \in \Hom(A^{n_1}, 
A^{m_1}),\: \ldots,\: \cA_i\in \Hom(A^{n_i},A^{m_i})$, the following diagram
commutes:
\[\xymatrix{
\rmcr_i(F)(A^{n_1}, \ldots, A^{n_i}) 
\ar[rrr]^-{\alpha_i(A^{n_1}, \ldots, A^{n_i})}
\ar[d]^-{\rmcr_i(F)(\cA_1, \ldots, \cA_i)} &&&
\rmcr_i(G)(A^{n_1}, \ldots, A^{n_i}) 
\ar[d]^-{\rmcr_i(G)(\cA_1, \ldots, \cA_i)} \\
\rmcr_i(F)(A^{m_1}, \ldots, A^{m_i})
\ar[rrr]^-{\alpha_i(A^{m_1}, \ldots, A^{m_i}) }&&&
\rmcr_i(G)(A^{m_1}, \ldots, A^{m_i}). }\]
For $i=d+1$, this assertion is clear since $\rmcr_{d+1}(F)$ and $\rmcr_{d+1}(G)$ are 
identically zero. So, we may assume that $i\le d$ and that the assertion is
already proved for $i+1$. Writing $(\cA_1, \ldots, \cA_i)$ as the composition
of $(\cA_1, \id, \ldots, \id),\: \ldots,\: (\id, \ldots, \id, \cA_i)$ and using
symmetry, we may furthermore assume that $\cA_2, \ldots, \cA_i$ are the
identity matrices. By construction of $\alpha_i(A^{n_1}, \ldots, A^{n_i})$, we
may finally assume that $n_2= \ldots = n_i =1$. Now, we write $\cA$ for
$\cA_1$, $\rmcr_i(F)(A^{n})$ for $\rmcr_i(F)(A^{n_1}, A, \ldots, A)$, and
similarly for $\rmcr_i(G)(A^{n_1}, A, \ldots, A)$, etc. By construction, the
above diagram commutes, if $\cA$ is a standard projection $A^n \ra A$ or a
standard embedding $A \ra A^m$. Furthermore, the diagram commutes if $n=m=1$
by assumption. Hence, the diagram commutes for all matrices $\cA$ which have
at most one entry which is different from zero. Thus it suffices to show
that the diagram commutes for the matrix $\cA + \cB$ if it commutes
for the matrices $\cA$ and $\cB$. For this, we decompose $\cA + \cB$ into the
composition 
\[\xymatrix@1{
A^n \ar[r]^-{\gD} & A^n \oplus A^n \ar[r]^-{\cA \oplus \cB} &
A^m \oplus A^m \ar[r]^-{+} & A^m }\]
and identify $\rmcr_i(-)(A^n \oplus A^n)$ with $\rmcr_i(-)(A^n) \oplus \rmcr_i(A^n)
\oplus \rmcr_{i+1}(-)(A^n, A^n)$ according to Proposition 1.2. By construction and
assumption, the following diagram commutes:
\[\xymatrix{
\rmcr_i(F)(A^n) 
\ar[rrrrr]^-{\alpha_i{(A^n)}}
\ar[d]^-{\rmcr_i(F)(\gD)} &&&&&
\rmcr_i(G)(A^n)
\ar[d]^-{\rmcr_i(G)(\gD)} \\
\rmcr_i(F)(A^n \oplus A^n)
\ar[rrrrr]^-{\alpha_i(A^n) \oplus \alpha_i(A^n) \oplus
\alpha_{i+1}(A^n,A^n)} &&&&&
\rmcr_i(G)(A^n \oplus A^n). }\]
(To see this, one decomposes the map $\rmcr_i(F)(\gD)$ into a direct sum of copies
of the maps
\[\rmcr_{i+j}(F)(\gD, \ldots, \gD): \rmcr_{i+j}(F)(A, \ldots, A) \ra
\rmcr_{i+j}(F)(A\oplus A, \ldots, A\oplus A), \quad j\ge 0, \]
according to Proposition 1.2. The composition of such a map with a natural projection
$\rmcr_{i+j}(F)(A\oplus A, \ldots, A\oplus A) \ra \rmcr_{i+j+k}(F)(A, \ldots, A)$
(again according to Proposition 1.2) is then a diagonal map $\gD_\varepsilon$ for
some $\varepsilon \in \{1,2\}^{i+j}$ with $|\varepsilon| = i+j+k$.) Likewise, the
corresponding diagram for $+$ in place of $\gD$ commutes. By the assumption
on $\cA$ and $\cB$ and by the induction hypothesis, also the following diagram
commutes:
\[\xymatrix{
\rmcr_i(F)(A^n \oplus A^n)
\ar[rrrrr]^-{\alpha_i(A^n) \oplus \alpha_i(A^n) \oplus \alpha_{i+1}(A^n, A^n)}
\ar[d]^-{\rmcr_i(F)(\cA \oplus \cB)} &&&&&
\rmcr_i(G)(A^n \oplus A^n)
\ar[d]^-{\rmcr_i(G)(\cA \oplus \cB)} \\
\rmcr_i(F)(A^m\oplus A^m) 
\ar[rrrrr]^-{\alpha_i(A^m) \oplus \alpha_i(A^m) \oplus \alpha_{i+1}(A^m, A^m)}&&&&&
\rmcr_i(G)(A^m\oplus A^m). }\]
Now, the proof of the above assertion is complete. \\
This assertion implies in particular that, for all $n\ge 1$ and for all 
projectors $\cA \in \End(A^n)$, the isomorphism
\[\xymatrix@1{
F(A^n) = \rmcr_1(F)(A^n) \ar[rr]^-{\alpha_1(A^n)} && 
\rmcr_1(G)(A^n) = G(A^n) }\]
induces an isomorphism
\[\xymatrix@1{
\alpha_1(\Image(\cA)):  F(\Image(\cA)) \cong \Image(F(\cA))
\ar[r]^-{\sim} & \Image(G(\cA)) \cong G(\Image(\cA)). }\] 
If $\cB \in \End{(A^m)}$ is a further projector and $f: \Image(\cA) \ra
\Image(\cB)$ is an $A$-homomorphism, then the diagram 
\[\xymatrix{
F(\Image(\cA)) 
\ar[rrr]^-{\alpha_1(\Image(\cA))} 
\ar[d]^-{F(f)} &&&
G(\Image(\cA)) 
\ar[d]^-{G(f)}\\
F(\Image(\cB))
\ar[rrr]^-{\alpha_1(\Image(\cB))} &&&
G(\Image(\cB)) }\]
commutes (again by the above assertion). Finally, if $V$ is an arbitrary
f.\ g.\ projective $A$-module, we choose a projector $\cA$ as above and
an isomorphism $V\cong \Image(\cA)$ and define $\alpha_1(V)$ to be
the composition
\[\xymatrix@1{
F(V) \cong F(\Image(\cA)) 
\ar[rrr]^-{\alpha_1(\Image(\cA))} &&&
G(\Image(\cA)) \cong G(V). }\]
It is clear that $\alpha_1(V)$ does not depend on the chosen projector and on the
chosen isomorphism and that $\alpha_1$ is an isomorphism of functors. This
proves Theorem 1.5.

{\bf Corollary 1.6} (Characterization of the exterior power functors). Let $A$
be a commutative ring, $B$ a further ring, $d \in \NN$, and 
\[ F : \cP_A := (\mbox{f.\ g.\ projective } A\mbox{-modules}) \ra
(B\mbox{-modules}) =: \cM_B\]
a functor with $F(0) =0$ and with the following properties:\\
(i) $\rmcr_d(F)$ is $d$-additive, i.e., $\rmcr_{d+1}(F)\equiv 0$. \\
(ii) The $d$ $A$-module structures on $\rmcr_d(F)(A, \ldots, A)$ coincide.\\
(iii) $F(A^{d-1}) =0$.\\
Then we have for all $V \in \cP_A$:
\[F(V) \cong \gL^d(V) \otimes_A \rmcr_d(F)(A, \ldots, A).\]

{\bf Proof.} The assumption $F(A^{d-1}) =0$ implies 
\[\rmcr_1(F)(A) =0, \:\ldots,\: 
\rmcr_{d-1}(F)(A, \ldots, A) = 0\] 
by Proposition 1.2. Similarly, for the functor
$G$ defined by $G(V) := \gL^d(V) \otimes_A \rmcr_d(V)(A, \ldots, A)$, we have
$\rmcr_1(G)(A) = 0, \:\ldots,\: \rmcr_{d-1}(G)(A, \ldots, A) =0$. Furthermore,
we have $\rmcr_d(G)(A, \ldots, A) \cong \rmcr_d(F)(A, \ldots, A)$. Thus, Corollary
1.6 follows from Theorem 1.5.

{\bf Example 1.7.} In the case $d =1$, Corollary 1.6 yields the following statement
which has already been proved in the papers \cite{E} and \cite{Wa}: Each additive
functor $F: \cP_A \ra \cM_B$ is isomorphic to the tensor functor $- \otimes F(A)$. 
By induction, we obtain the following generalization: Let $F: \cP_A^d \ra \cM$
be a $d$-additive functor such that the $d$ $A$-module structures on 
$F(A, \ldots, A)$ coincide. Then we have for all $V_1, \ldots, V_d \in \cP_A$:
\[ F(V_1, \ldots, V_d) \cong V_1 \otimes_A \ldots \otimes_A V_d \otimes_A
F(A, \ldots, A).\]

{\bf Proposition 1.8} (Characterization of the symmetric power functors). Let 
$A$, $B$, $d$, and $F$ as in Corollary 1.6. In contrast to Corollary 1.6, we assume 
however that $F$ has the following property
in place of property (iii): For all $i = 1, \ldots, d-1$, the $B$-module
homomorphism
\[\xymatrix{
{\rm plus}_i: \oplusm_{\varepsilon \in \{1, \ldots, d\}^i , \, |\varepsilon| =d}
\rmcr_d(F)(A, \ldots, A) 
\ar[rrrr]^-{(+_\varepsilon)_{\varepsilon \in \{1, \ldots, d\}^i, 
|\varepsilon|=d }} &&&&
\rmcr_i(F)(A, \ldots, A) }\]
is bijective. Then we have for all $V \in \cP_A$:
\[ F(V) \cong \Sym^d(V) \otimes_A \rmcr_d(F)(A, \ldots, A).\]

{\bf Proof.} We view $\rmcr_d(F)(V, \ldots, V)$ as a $\gS_d$-module with the
action introduced above and $F(V)$ as a $\gS_d$-module with trivial action. Then,
the plus map 
\[+_d: \rmcr_d(F)(V, \ldots, V) \ra F(V)\] 
is obviously compatible with
these $\gS_d$-actions. Since $+_d: \rmcr_d(F)(A, \ldots, A) \ra F(A)$ is 
bijective by assumption, the symmetric group $\gS_d$ acts trivially on 
$\rmcr_d(F)(A, \ldots, A)$. Hence, the composition 
\[V^{\otimes d} \otimes_A \rmcr_d(F)(A, \ldots, A) \cong \rmcr_d(F)(V, \ldots, V)
\,\, \stackrel{+_d}{\longrightarrow} \,\,  F(V) \]
of the isomorphism given in Example 1.7 with the plus map $+_d$ induces a 
$B$-module homomorphism 
\[\alpha(V) : \Sym^d(V) \otimes_A \rmcr_d(F)(A, \ldots, A) \ra F(V)\]
which is functorial in $V$. One easily proves that, for all $i= 1, \ldots, d$, 
the homomorphism $\rmcr_i(\alpha)(A, \ldots, A)$ coincides with the isomorphism 
${\rm plus}_i$ by virtue of the isomorphism 
\[\rmcr_i(\Sym^d)(A, \ldots, A) \cong 
\oplusm_{\varepsilon \in \{1, \ldots, d\}^i, \, |\varepsilon| =d} A\]
(see Example 1.4(a)). Hence, $\rmcr_i(\alpha)(A, \ldots, A)$ is bijective.
Then, by Proposition 1.2, $\alpha(V)$ is bijective for all f.\ g.\ {\em free}
$A$-modules $V$. Again by Proposition 1.2, the same then holds for {\em all}
$V \in \cP_A$. This proves Proposition 1.8.

{\bf Proposition 1.9} (Characterization of the divided power functors).
Let $A$, $B$, $d$, and $F$ be as in Corollary 1.6. In contrast to Corollary 1.6,
we assume however that $F$ has
the following property in place of property (iii): For all $i=1, \ldots, d-1$,
the $B$-module homomorphism
\[\xymatrix@1{
{\rm diag}_i: \rmcr_i(F)(A, \ldots, A)  
\ar[rrrr]^-{(\gD_\varepsilon)_{\varepsilon\in \{1, \ldots, d\}^i, 
|\varepsilon|=d}} &&&& \oplusm_{\varepsilon\in \{1, \ldots, d\}^i, 
|\varepsilon|=d} \rmcr_d(F)(A, \ldots, A)}\]
is bijective. Then we have for all $V \in \cP_A$:
\[F(V) \cong D^d(V) \otimes_A \rmcr_d(F)(A, \ldots, A).\]

{\bf Proof.} The $d$-th divided power $D^d(V)$ of a f.\ g.\ projective $A$-module
$V$ is the module of fixed elements under the natural action of the symmetric
group $\gS_d$ on $V^{\otimes d}$. Thus, Proposition 1.9 can be proved analogously
to Proposition 1.8.

{\bf Remark 1.10.} For any functor $F: \cP_A \ra \cM_B$ with $F(0) =0$, the
composition $\xymatrix@1{F(A) \ar[r]^-{\gD_2} &
\rmcr_2(F)(A,A) \ar[r]^-{+_2} & F(A) }$ equals $F(2 \cdot \id_A) - 2 \cdot 
\id_{F(A)}$. Now assume that $d=2$. 
Under the assumptions of Proposition 1.8, the endomorphism $F(2 \cdot \id_A)$
of $F(A)$ then corresponds to the endomorphism $\rmcr_2(F)(2 \cdot \id_A, 
2\cdot \id_A) = 4 \cdot \id$ of $\rmcr_2(F)(A,A)$; hence, $F(2 \cdot \id_A)=
4 \cdot \id_{F(A)}$. Thus, $+_2 \circ \gD_2 = 2 \cdot \id_{F(A)}$. This proves that
the assumptions of Proposition 1.8 also determine the map $\gD_2$. Thus, 
Proposition 1.8 follows from Theorem 1.5 in the case $d=2$. The same holds for
Proposition 1.9. It seems that a simple generalization of this argument to the 
case $d \ge 2$ does not exist. 

\bigskip

\section*{\S 2 Simplicial Modules and Koszul Complexes}

First, we recall the Dold-Kan correspondence; in particular, we introduce the
functors $\gG$ and $N$ between the category of simpicial modules and the category
of complexes. Then, we recall the definition of the $n$-th Koszul complex
$\Kos^n(f)$ associated with a homomorphism $f: P \ra Q$ between the projective
modules $P,Q$ over a commutative ring $R$. Finally, we construct an explicit
complex homomorphism $u^n(f)$ from $\Kos^n(f)$ to $N \Sym^n \gG(\ldots \ra 
P \,\, \stackrel{f}{\ra}\,\, Q)$ and prove that $u^n(f)$ is a quasi-isomorphism.

Let $\cM$ be an additive category. Let $\gD$ denote the category whose objects are
the sets $[n]=\{0 < 1 < \ldots < n\}$, $n\ge 0$, and whose morphisms are the 
order preserving set maps. For $0\le i \le n+1$, let $\gd_i^n: [n] \ra [n+1]$
denote the injective order preserving map given by $\Image(\gd_i^n) =
[n+1] \backslash \{i\}$. For $0\le i \le n$, let $\gs_i^n: [n+1] \ra [n]$ denote
the surjective order preserving map given by $(\gs_i^n)^{-1}(i) = \{i, i+1\}$. 
A {\em simplicial object $X.$ in $\cM$} is a contravariant functor from 
$\gD$ to $\cM$. We write $X_n$ for $X.([n])$, $d_i: X_{n+1} \ra X_n$ for 
$X.(\gd_i^n)$ and $s_i: X_n \ra X_{n+1}$ for $X.(\gs_i^n)$. A {\em simplicial
homomorphism} $f. : X. \ra Y.$ between two simplicial objects $X., Y.$ in 
$\cM$ is a natural transformation. A {\em simplicial homotopy} between 
simplicial homomorphisms $f,g: \xymatrix{X. \ar@<0.5ex>[r] \ar@<-0.5ex>[r] & Y}.$
consists of homomorphisms $h_i(n): X_n \ra Y_{n+1}$, $n\ge 0$, $0\le i \le n$, 
such that $d_0 h_0(n) =f_n$ and $d_{n+1} h_n(n) = g_n$ for all $n\ge 0$ and such
that certain further relations $d_i h_j(n) = \ldots$ and $s_ih_j(n) =
\ldots$ hold (see 8.3.11 in \cite{We} or (2.3) in \cite{JM}). A {\em complex 
$K.$ in $\cM$} is a contravariant functor from the ordered set $\NN_0$ to
the category $\cM$. The {\em normalized complex} $N(X.)$ of a simplicial
object $X.$ in an abelian category $\cM$ is given by 
\[N(X.)_n := X_n \left/ \sum_{i=0}^{n-1} \Image(s_i: X_{n-1} \ra X_n) \right.;\]
the differential in $N(X.)$ is 
\[\partial = \sum_{i=0}^{n} (-1)^{i} d_i : N(X.)_{n} \ra N(X.)_{n-1}\]
(for all $n\ge 0$). Thus, we have a functor
\[N: \Simp(\cM) \ra \Compl(\cM)\]
from the category $\Simp(\cM)$ of simplicial objects in $\cM$ to the category
$\Compl(\cM)$ of complexes in $\cM$. We define a functor
\[ \gG: \Compl(\cM) \ra \Simp(\cM)\]
in the reverse direction as follows. For $n\ge k$, let $\Sur([n],[k])$ denote the
set of surjective order preserving maps from $[n]$ to $[k]$. For $(M., \partial)
\in \Compl(\cM)$ and $n\ge 0$, let $\gG(M.)_n$ be defined by
\[\gG(M.)_n := \oplusm_{k=0}^n \oplusm_{\gs \in \Sur([n],[k])} M_k;\]
furthermore, $d_i: \gG(M.)_n \ra \gG(M.)_{n-1}$, $i \in \{0, \ldots, n\}$, and
$s_i: \gG(M.)_n \ra \gG(M.)_{n+1}$, $i \in \{0, \ldots, n\}$, are defined by
\begin{eqnarray*}
d_i(m,\gs) & := & 
\left\{\begin{array}{ll}
(m,\gs\gd_i) & \textrm{if } \gs\gd_i \in \Sur([n-1],[k])\\
(\partial(m), \hat{\gs}) & \textrm{if } \gs\gd_i = \gd_0\hat{\gs} \: 
(\textrm{where } \hat{\gs} \in \Sur([n-1],[k-1]))\\
0 & \textrm{if } \gs\gd_i = \gd_j \hat{\gs} \textrm{ with } j\ge 1
\end{array}\right.\\
s_i(m, \gs) & := & (m, \gs\gs_i)
\end{eqnarray*}
(for $k\in \{0, \ldots, n\}$ and $(m,\gs) \in M_k \times \Sur([n],[k])$). One 
easily verifies that these definitions determine a simplicial object 
$\gG(M.) $ in $\cM$ (cf.\ section 2 in \cite{JM} or section 8.4 in \cite{We}). 

{\bf Proposition 2.1} (Dold-Kan correspondence). For any abelian category $\cM$,
the functors $N$ and $\gG$ are inverse to each other (up to canonical isomorphisms). 
A simplicial homotopy between two simplicial homomorphisms induces a homotopy
between the corresponding complex homomorphisms and vice versa.

{\bf Proof.} See section 8.4 in \cite{We}.

Now, let $F: \cP \ra \cM$ be a functor from an additive category $\cP$ to an
abelian category $\cM$ with $F(0) = 0$. If $F$ is additive, then the induced
functor $\Compl(\cP) \ra \Compl(\cM)$, $K. \mapsto F\circ K.$, maps homotopic
complex homomorphisms to homotopic complex homomorphisms. In the simplicial world,
the analogous statement holds even for an arbitrary functor $F$: The induced functor
$\Simp(\cP) \ra \Simp(\cM)$, $X. \mapsto FX. := F\circ X.$, maps simplicial
homotopies to simplicial homotopies since only compositions (and no sums) are involved
in the defining relations of a simplicial homotopy. Hence, using the Dold-Kan
correspondence, we obtain a functor
\[\Compl(\cP) \ra \Compl(\cM), \quad K. \mapsto NF\gG(K.),\]
which maps homotopic complex homomorphisms to homotopic complex homomorphisms and
which is obviously isomorphic to the functor $K. \mapsto F \circ K.$, if 
$F$ is additive. (This basic fact has already been observed by Dold and Puppe 
in \cite{DP}.) If $K.$ is a complex in $\cP$ with $K_n = 0$ for all $n\ge 2$,
then the complex $NF\gG(K.)$ can be described as follows by using the 
cross effect functors introduced in section 1.

{\bf Lemma 2.2.} Let $f: P\ra Q$ be a homomorphism in $\cP$. We consider $f$ as
a complex with $Q$ at the place $0$ and $P$ at the place $1$. Then we have for
all $n\ge 0$:
\[NF\gG(P\ra Q)_n \cong \rmcr_n(F)(P,\ldots, P) \oplus \rmcr_{n+1}(F)(Q, P, \ldots, P);\]
the differential $\partial: NF\gG(P\ra Q)_n \ra NF\gG(P\ra Q)_{n-1}$ is given
by the following diagram:
\[\xymatrix{
\rmcr_n(F)(P, \ldots, P) \ar[r] \ar[rd] \ar@{}[d]|{\oplus} &
\rmcr_{n-1}(F)(P,\ldots, P) \ar@{}[d]|{\oplus} \\
\rmcr_{n+1}(F)(Q, P, \ldots, P) \ar[r] & \rmcr_n(F)(Q,P,\ldots, P);
}\]
here, the upper horizontal map is the alternating sum $\sum_{i=1}^{n-1}
(-1)^i (+_{\varepsilon^i})$ of the plus maps associated with 
$\varepsilon^1= (2,1,\ldots, 1), \: \ldots, \: \varepsilon^{n-1} = (1, \ldots, 1,2)$;
the lower horizontal map is $\left((+_{(2,1,\ldots,1)}) \circ \rmcr_{n+1}(F)
(\id_Q, f, \id_P, \ldots, \id_P) \right) + \sum_{i=1}^{n-1}(-1)^i(+_{(1,\varepsilon^i)})$;
the diagonal map is $\rmcr_n(F)(f, \id_P, \ldots, \id_p)$. 

{\bf Proof.} By construction of $\gG$, we have $\gG(P\ra Q)_n \cong Q \oplus P^n$
for all $n\ge 0$. By virtue of these isomorphisms, the maps $d_i:
\gG(P\ra Q)_n \ra \gG(P\ra Q)_{n-1}$ and $s_i:\gG(P\ra Q)_n \ra 
\gG(P\ra Q)_{n+1}$ can be described as follows:
\begin{eqnarray*}
d_i(q,p_1, \ldots, p_n) & = & \left\{
\begin{array}{ll}
(q+f(p_1), p_2, \ldots, p_n) & \textrm{for } i=0\\
(q, p_1, \ldots, p_{i-1}, p_i+p_{i+1}, p_{i+2}, \ldots, p_n) &
\textrm{for } 1 \le i \le n-1\\
(q,p_1, \ldots, p_{n-1}) & \textrm{for } i=n
\end{array}\right. \\
s_i(q, p_1, \ldots, p_n) &= &(q, p_1, \ldots, p_i, 0, p_{i+1}, \ldots, p_n) \quad
\textrm{for } 0\le i \le n
\end{eqnarray*}
(see also Lemma A.2 in \cite{JM}). Hence, by Proposition 1.2, we have:
\begin{eqnarray*}
\lefteqn{F\gG(P\ra Q)_n \cong F(Q) \oplus F(P^n) \oplus \rmcr_2(F)(Q, P^n) }\\
&& \cong \oplusm_{k=0}^n \oplusm^{{n \choose k}}
\left( \rmcr_k(F)(P, \ldots, P) \oplus \rmcr_{k+1}(F)(Q, P, \ldots, P)\right).
\end{eqnarray*}
Thus we have
\begin{eqnarray*}
\lefteqn{NF\gG(P \ra Q)_n \cong F\gG(P \ra Q)_n 
\left/ \left( \sum_{i=0}^{n-1}\Image(s_i)\right)\right.} \hspace{2cm}\\
&& \cong \rmcr_n(F)(P, \ldots, P) \oplus \rmcr_{n+1}(F)(Q,P, \ldots, P).
\end{eqnarray*}
Finally, one can easily check now that the differential has the claimed form.

For any object $P$ in $\cP$ and $k\ge 0$, we write $P[-k]$ for the complex
which has $P$ at the place $k$ and $0$ else. Lemma 2.2 states in particular
that $NF\gG(P[0]) \cong F(P)[0]$ and that $NF\gG(P[-1])$ is isomorphic to
the complex
\[\ldots \longrightarrow \rmcr_n(F)(P, \ldots, P) \,\, 
\stackrel{\partial_n}{\longrightarrow}\,\,
\rmcr_{n-1}(F)(P, \ldots, P) \,\, \stackrel{\partial_{n-1}}{\longrightarrow} \,\,
\ldots \,\, \stackrel{\partial_2}{\longrightarrow} \,\, F(P) \longrightarrow 0\]
where $\partial_n= \sum_{i=1}^{n-1}(-1)^i(+_{\varepsilon^i})$.

Now, we recall the definition of the Koszul complex (see I 4.3.1.3 in \cite{Il},
V.1.3 in \cite{ABW}, or  section 2 in \cite{Gr}). Let $n\in \NN$, $R$ a
commutative ring, and $f:P\ra Q$ a homomorphism between projective $R$-modules $P,Q$.

{\bf Definition 2.3.} For any $k \in \{0, \ldots, n+1\}$, let 
\[d_{k+1}: \gL^{k+1}(P) \otimes \Sym^{n-k-1}(Q) \ra 
\gL^k(P) \otimes \Sym^{n-k}(Q) \]
denote the Koszul differential given by
\[p_1 \wedge \ldots \wedge p_{k+1} \otimes q_{k+2} \cdots q_n 
\mapsto  \sum_{i=1}^{k+1} (-1)^{k+1-i} p_1 \wedge \ldots \wedge \hat{p}_i \wedge
\ldots \wedge p_{k+1} \otimes f(p_i) q_{k+2}\cdots q_n.\]
The complex 
\[ 0 \ra \gL^n(P) \,\, \stackrel{d_n}{\ra} \,\, 
\gL^{n-1}(P)\otimes Q \,\, \stackrel{d_{n-1}}{\ra} \,\, \ldots \,\,
\stackrel{d_2}{\ra} P \otimes \Sym^{n-1}(Q) \,\, 
\stackrel{d_1}{\ra}\,\, \Sym^n(Q) \ra 0\]
is called {\em $n$-th Koszul complex associated with $f$} and is denoted by
$\Kos^n(f)$.

One can prove analogously to Lemme 2.1.2.1 on page 277 in \cite{Il2} that the
complexes $\Kos^n(f)$ and $N \Sym^n \gG(P \ra Q)$ are isomorphic in the 
derived category of $R$. Whereas Illusie uses three complex homomorphisms (which
even can not be composed since they have different directions) to realize this
isomorphism, we now explicitly construct a single complex homomorphism $u^n(f)$
from $\Kos^n(f)$ to $N \Sym^n \gG(P\ra Q)$ and show that $u^n(f)$ is a 
quasi-isomorphism. For $k=0, \ldots, n$, let $u_k^n(f)$ denote the composition
\begin{eqnarray*}
\lefteqn{\gL^k(P) \otimes \Sym^{n-k}(Q) 
\,\, \stackrel{a_k \otimes \id}{\longrightarrow} \,\,
P^{\otimes k} \otimes \Sym^{n-k}(Q) }\\
&& \stackrel{(1.4)(a)}{\tilde{\longrightarrow}} \,\,
\rmcr_k(\Sym^k)(P, \ldots, P) \otimes \Sym^{n-k}(Q) \\
&& \stackrel{{\rm can}}{\hookrightarrow} \,\, 
\oplusm_{i=0}^{n-k} \rmcr_k(\Sym^{k+i})(P, \ldots, P) \otimes \Sym^{n-k-i}(Q)\\
&& \stackrel{(1.4)(a)}{\tilde{\longrightarrow}}\,\, 
\rmcr_k(\Sym^n)(P, \ldots, P) \oplus \rmcr_{k+1}(\Sym^n)(Q,P,\ldots, P) \\
&& \stackrel{(2.2)}{\tilde{\longrightarrow}}\,\, N \Sym^n\gG(P\ra Q)_k;
\end{eqnarray*}
here, $a_k : \gL^k(P) \ra P^{\otimes k}$, $p_1\wedge \ldots \wedge p_k \mapsto
\sum_{\gs \in \gS_k} (-1)^{\sgn(\gs)} p_{\gs(1)} \otimes \ldots \otimes 
p_{\gs(k)},$ denotes the antisymmetrization map; the other homomorphisms are
explained under the indicated example or lemma.

{\bf Proposition 2.4.} The homomorphisms $u_k^n(f)$, $k=0, \ldots, n$, define a 
quasi-iso\-morphism
\[u^n(f): \Kos^n(f) \,\, \tilde{\ra} \,\, N\Sym^n\gG(P\ra Q).\]

{\bf Proof.} We set $L. := N\Sym^n\gG(P\ra Q)$. The following picture 
illustrates $u^n(f)$:
\[\xymatrix{
\gL^n(P) \ar[r] \ar[ddddd]^-{a_n} & 
\gL^{n-1}(P) \otimes Q \ar[r] \ar[dddd]^-{a_{n-1}\otimes \id} &
\ldots \ar[r] & P\otimes \Sym^{n-1}(Q) \ar[r] \ar[dd]^-{a_1 \otimes \id} & 
\Sym^n(Q) \ar[d]^-{\id}\\
&&&& \Sym^n(Q) \ar@{}[d]|{\oplus} \\
&&& P\otimes \Sym^{n-1}(Q) \ar[ur] \ar[r] \ar@{}[d]|{\oplus} &
0 \ar@{}[d]|{\oplus} \\
&& \antiddots  {}\ar[ur]& \vdots \ar@{}[d]|{\oplus} &  \vdots \ar@{}[d]|{\oplus}\\
& P^{\otimes (n-1)} \otimes Q \ar[ur] \ar[r] \ar@{}[d]|{\oplus} &
\ldots \ar[ur] \ar[r] & \Sym^{n-1}(P)\otimes Q \ar[r] \ar@{}[d]|{\oplus} & 
0 \ar@{}[d]|{\oplus} \\
P^{\otimes n} \ar[ur] \ar[r] \ar@{=}[d] &
\rmcr_{n-1}(\Sym^n)(P, \ldots, P) \ar[ur] \ar[r] \ar@{=}[d] &
\ldots \ar[ur] \ar[r] & \Sym^n(P)  \ar[r] \ar@{=}[d] & 0 \ar@{=}[d] \\
L_n \ar[r]^-{\partial} &
L_{n-1} \ar[r]^-{\partial} &
\ldots \ar[r]^-{\partial} & L_1 \ar[r]^-{\partial} &
L_0;}\]
here, the differentials in the middle complex are given by Lemma 2.2.
For all $k=0, \ldots, n$, the composition
\begin{eqnarray*}
\lefteqn{\xymatrix{
\gL^k(P) \ar[r]^-{a_k}& P^{\otimes k} \cong \rmcr_k(\Sym^k)(P, \ldots , P)
\ar[rrr]^-{\sum_{i=1}^{k-1} (-1)^i(+_{\varepsilon^i})} &&&
\rmcr_{k-1}(\Sym^k)(P, \ldots, P) }} \hspace{2cm} \\
&& \cong \left(\Sym^2(P) \otimes P \otimes \ldots \otimes P\right) \oplus
\ldots \oplus \left( P \otimes \ldots \otimes P \otimes \Sym^2(P)\right)
\end{eqnarray*}
is obviously the zero map; furthermore, one easily sees that the diagram
\[\xymatrix{
\gL^{k+1}(P) \otimes \Sym^{n-k-1}(Q) \ar[rr]^-{d_{k+1}} 
\ar[d]^-{a_{k+1} \otimes \id} &&
\gL^k(P) \otimes \Sym^{n-k}(Q) \ar[d]^-{a_k \otimes \id}\\
P^{\otimes (k+1)} \otimes \Sym^{n-k-1}(Q) 
\ar[r]^-{\id \otimes f \otimes \id} & 
P^{\otimes k} \otimes Q \otimes \Sym^{n-k-1}(Q) 
\ar[r]^-{{\rm can}} &
P^{\otimes k} \otimes \Sym^{n-k}(Q)
}\]
commutes. This proves that $u$ is a complex homomorphism. Obviously, the complexes 
$\Kos^n(f)$ and $N\Sym^n\gG(P\ra Q)$ can be filtered by subcomplexes such that
$u^n(f)$ is compatible with these filtrations and such that, for $k=0, \ldots, n$,
the map between the $k$-th successive quotients induced by $u^n(f)$ is the complex
homomorphism
\[u^k(0) \otimes \id : \gL^k(P)[-k] \otimes \Sym^{n-k}(Q) \ra
N\Sym^k\gG(P\ra 0) \otimes \Sym^{n-k}(Q).\]
Thus it suffices to show that, for all $n\ge 0$, the complex homomorphism
\[u^n(0): \gL^n(P)[-n] \ra N\Sym^n\gG(P\ra 0)\]
is a quasi-isomorphism, i.e., that Proposition 2.4 is true in the case $Q=0$. 
We prove this by induction on $n$. For $n=0$ this is clear. So let $n\ge 1$.
The homomorphism $u^n(\id_P)$ is a quasi-isomorphism since 
$\Kos^n(\id_P)$ is exact and the complex $N\Sym^n\gG(P\,\,\stackrel{\id}{\ra}\,\,P)$
is homotopic to the zero complex. Furthermore, the complex homomorphisms
\[u^k(0) \otimes \id: \gL^k(P)[-k] \otimes \Sym^{n-k}(P) \ra
N\Sym^k\gG(P\ra 0) \otimes \Sym^{n-k}(P), \]
$k=0, \ldots, n-1,$
are quasi-isomorphisms by the induction hypothesis. Using the filtrations already
considered above, we obtain that $u^n(0)$ is a quasi-isomorphism. (Alternatively,
this also follows from Proposition I 4.3.2.1(i) on p.\ 111 in \cite{Il} and the
following fact which is a little bit tedious to prove: 
The complex homomorphism $N\Sym^n\gG(P\ra 0) \ra \gL^n(P)[-n]$ 
constructed in \cite{Il} is left inverse to $u^n(0)$.) This ends the proof of 
Proposition 2.4.

Any commutative square
\[\xymatrix{
P \ar[r]^-{f} \ar[d]^-{\alpha_P} & Q \ar[d]^-{\alpha_Q}\\
P' \ar[r]^-{f'} & Q'
}\]
of homomorphisms between projective $R$-modules induces a homomorphism
\[\Kos^n(\alpha): \Kos^n(f) \ra \Kos^n(f')\]
between the corresponding Koszul complexes in the obvious way.

{\bf Corollary 2.5.} If $(\alpha_P^1, \alpha_Q^1)$ and $(\alpha_P^2, \alpha_Q^2)$
are homotopic complex homomorphisms between the complexes $P\,\, 
\stackrel{f}{\ra} \,\, Q$ and $P' \,\, \stackrel{f'}{\ra}\,\, Q'$, then
the induced homomorphisms between the $k$-th homology modules of the 
corresponding $n$-th Koszul complexes are equal:
\[H_k(\Kos^n(\alpha^1)) = H_k(\Kos^n(\alpha^2)):
H_k(\Kos^n(f)) \ra H_k(\Kos^n(f'))\]
(for all $k\ge 0$).

{\bf Proof.} Obviously, the diagram
\[\xymatrix{
\Kos^n(f) \ar[rr]^-{u^n(f)} \ar[d]^-{\Kos^n(\alpha^i)} &&
N\Sym^n\gG(P\,\, \stackrel{f}{\ra}\,\, Q) \ar[d]^-{N\Sym^n\gG(\alpha^i)}\\
\Kos^n(f') \ar[rr]^-{u^n(f')} &&
N\Sym^n\gG (P' \,\,\stackrel{f'}{\ra}\,\, Q')
}\]
commutes for $i=1,2$. Furthermore, the complex homomorphisms 
$N\Sym^n\gG(\alpha^1)$ and $N\Sym^n\gG(\alpha^2)$ are homotopic to each other.
Thus, Corollary 2.5 follows from Proposition 2.4.

{\bf Remark 2.6.} In the next section, we will use only Corollary 2.5 (and not
Proposition 2.4). The simplicial methods used here to prove Corollary 2.5 are actually
not necessary: Alternatively, one can use use Lemma 2.2 to define
the complex $NF\gG(P\ra Q)$ and one can directly show that homotopic complex 
homomorphisms between $P\ra Q$ and $P'\ra Q'$ induce homotopic complex 
homomorphisms between $NF\gG(P\ra Q)$ and $NF\gG(P'\ra Q')$. However, the simplicial
definition of $NF\gG(K.)$ has the advantage that it yields a more natural
argument for the proof of Corollary 2.5 and that it can be applied also to
complexes $K.$ of length $>1$ (see also section 6).

{\bf Remark 2.7.} One can prove analogously to Proposition 2.4 that the Koszul complex
\[\tilde{\Kos}^n(f): 0 \ra D^n(P) \ra D^{n-1}(P) \otimes Q \ra \ldots \ra 
P\otimes \gL^{n-1} (Q) \ra \gL^n(Q) \ra 0\]
(defined analogously to 2.3) is quasi-isomorphic to the complex
$N\gL^n\gG(P\ra Q)$.

\bigskip

\section*{\S 3 Computation of the Homology of Koszul Complexes}

Let $R$ be a commutative ring and $n\in \NN$. In this section, we compute the
homology modules of the $n$-th Koszul complex $\Kos^n(f)$ for each homomorphism
$f: P \ra Q$ between f.\ g.\ projective $R$-modules which has the following
properties: $f$ is injective and $Q/f(P)$ is a f.\ g.\ projective 
$R/I$-module for some ideal $I$ in $R$ which is locally generated by a 
non-zero-divisor.

First, we recall the definition of the Schur modules corresponding to 
Young diagrams of hook type.

{\bf Definition 3.1.} Let $k\in \{0, \ldots, n\}$ and $V$ a f.\ g.\ projective
$R$-module. Let 
\[d_{k+1}: \gL^{k+1}(V) \otimes \Sym^{n-k-1}(V) \ra 
\gL^k(V) \otimes \Sym^{n-k}(V)\] 
denote the differential in the Koszul complex
$\Kos^n(\id_V)$ (see Definition 2.3). The module
\[L_k^n(V) := \Image(d_{k+1}) \subseteq \gL^k(V) \otimes \Sym^{n-k}(V)\]
is called {\em Schur module associated with $V$ and the Young diagram
$(k+1, 1, \ldots, 1)$} ($(n-k-1)$ $1$s).

Since $\Kos^n(\id_V)$ is exact, $L^n_k(V)$ is a f.\ g.\ projective $R$-module
for all $k=0, \ldots, n$. For instance, we have $L_0^n(V) \cong \Sym^n(V)$, 
$L_{n-1}^n(V) \cong \gL^n(V)$, and $L_n^n(V) \cong 0$. Obviously, the definition
given here agrees with the definition on p.\ 220 in \cite{ABW}. 

{\bf Theorem 3.2.} Let $I$ be an ideal in $R$ which is locally generated by
a non-zero-divisor. Let $0\ra P \,\, \stackrel{f}{\ra}\,\, Q \ra V \ra 0$ be
an $R$-projective resolution of a f.\ g.\ projective $R/I$-module $V$. Then,
for all $k=0, \ldots, n$, we have a canonical isomorphism
\[H_k(\Kos^n(f)) \cong L^n_k(V) \otimes (I/I^2)^{\otimes k}\]
between the $k$-th homology module of the $n$-th Koszul complex and the 
tensor product of the Schur module $L^n_k(V)$ with the $k$-th tensor power of
the conormal module $I/I^2$. In particular, $H_k(\Kos^n(f))$ is a projective
$R/I$-module.

{\bf Remark 3.3.} We consider the following extreme cases:\\
(a) If $k=n$, Theorem 3.2 states that the differential
\[d_n: \gL^n(P) \ra \gL^{n-1}(P) \otimes Q\]
is injective. This immediately follows from the exactness of the Koszul complex
$\Kos^n(\id_P)$ and the injectivity of $f$.\\
(b) If $k=0$, Theorem 3.2 states that the sequence
\[P\otimes \Sym^{n-1}(Q) \,\, \stackrel{d_1}{\longrightarrow} \,\,
\Sym^n(Q) \,\, \stackrel{{\rm can}}{\longrightarrow} \,\,
\Sym^n(V) \longrightarrow 0\]
is exact. This already follows from the construction of symmetric powers.\\
In both cases, the assumption on the ideal $I$ is not used.

{\bf Proof of Theorem 3.2.} We fix $k\in \{0, \ldots, n\}$. For any f.\ g.\
projective $R/I$-module $V$, we choose an $R$-projective resolution 
$0 \ra P_V \,\, \stackrel{f_V}{\ra} Q_V \ra V \ra 0$ of $V$. For any 
homomorphism $\alpha: V\ra W$ between
f.\ g.\ projective $R/I$-modules $V,W$, there is a complex homomorphism
\[\xymatrix{
P_V \ar[r] \ar[d]^-{\alpha_P} & Q_V \ar[d]^-{\alpha_Q}\\
P_W \ar[r] & Q_W
}\]
which is compatible with $\alpha$; this complex homomorphism is unique up
to homotopy by a standard argument of homological algebra. Thus, by Corollary 2.5, we obtain a well-defined functor
\[F: \cP_{R/I} \ra (R\textrm{-modules}), \quad V \mapsto
H_k(\Kos^n(f_V)).\]
Let $G$ denote the functor
\[G: \cP_{R/I} \ra (R\textrm{-modules}), \quad V \mapsto
L^n_k(V) \otimes (I/I^2)^{\otimes k}.\]
We have to show that $F(V) \cong G(V)$ for all $V\in \cP_{R/I}$. For this, we
verify the assumptions of Theorem 1.5.\\
One easily checks (see also 4.3.1.5 on p.\ 109 in \cite{Il}) that, for all
$f:P\ra Q$, $f': P'\ra Q'$ in $\cP_R$, we have:
\[\Kos^n(f\oplus f') \cong \oplusm_{n_1 + n_2=n}
\Kos^{n_1}(f) \otimes \Kos^{n_2}(f').\]
(The tensor product of two complexes $K.$ and $L.$ is a priori only a
double complex. In this paper, $K. \otimes L.$ denotes both the double complex and the 
total complex of this double complex.) Using Proposition 1.2, we obtain the 
following isomorphism for all $i \ge 0$ and $V_1, \ldots, V_i \in \cP_{R/I}$
(see also Example 1.4(a)):
\[\rmcr_i(F)(V_1, \ldots, V_i) \cong \oplusm_{{(n_1, \ldots, n_i) \in \{1, \ldots, n\}^i
\atop n_1 + \ldots + n_i =n}} H_k\left(\Kos^{n_1}(f_{V_1}) \otimes \ldots \otimes 
\Kos^{n_i}(f_{V_i})\right).\]
In particular, $F$ is of degree $\le n$ and, for $V_1 = \ldots = V_i= R/I$, we obtain
\[\rmcr_i(F)(R/I, \ldots, R/I) \cong \oplusm_{{(n_1, \ldots, n_i) \in \{1, \ldots, n\}^i
\atop n_1 + \ldots + n_i =n}} H_k\left( (I \,\, \stackrel{{\rm can}}{\longrightarrow}\,\,
R)^{\otimes i} \right)\]
since $\Kos^m(I \ra R) \cong (I\ra R)$ for all $m\ge 1$. 
The complex $(I \,\,\stackrel{{\rm can}}{\longrightarrow} \,\, R) ^{\otimes i}$ has
the shape
\[0 \ra I^{\otimes i} \ra \oplusm^{{i \choose i-1}} I^{\otimes (i-1)} \ra \ldots
\ra \oplusm^{{i \choose 1}} I \ra R \ra 0;\]
here, the differential $\oplus^{{i \choose k}} I^{\otimes k} \ra 
\oplus^{{i \choose k-1}} I^{\otimes (k-1)}$ is the composition of the differential
$\oplus^{{i \choose k}} I^{\otimes k} \ra \oplus^{{i \choose k-1}} I^{\otimes k}$
in the complex $(R \, \, \stackrel{\id}{\ra} R)^{\otimes i} \otimes I^{\otimes k}$
with the canonical inclusion $\oplus^{{i \choose k-1}} I^{\otimes k} \ra 
\oplus^{{i \choose k-1}} I^{\otimes (k-1)}$. Hence, the module
$Z_k\left((I \ra R)^{\otimes i}\right)$ of $k$-cycles is isomorphic to
$Z_k\left( (R\ra R)^{\otimes i}\right) \otimes I^{\otimes k}$, and the module
$B_k\left((I\ra R)^{\otimes i}\right)$ of $k$-boundaries is isomorphic to
$B_k\left((R\ra R)^{\otimes i}\right) \otimes I^{\otimes (k+1)}$. Since 
$(R\ra R)^{\otimes i}$ is exact, we obtain:
\[H_k\left((I\ra R)^{\otimes i}\right) \cong B_k\left((R\ra R)^{\otimes i}\right)
\otimes (I/I^2)^{\otimes k}.\]
Similarly, for the functor $G$, we obtain:
\[\rmcr_i(G)(V_1, \ldots, V_i) \cong 
\oplusm_{{(n_1, \ldots, n_i) \in \{1, \ldots, n\}^i \atop n_1 + \ldots + n_i =n}}
B_k\left(\Kos^{n_1}(\id_{V_1}) \otimes \ldots \otimes \Kos^{n_i}(\id_{V_i})\right)
\otimes (I/I^2)^{\otimes k}.\]
In particular, $G$ is of degree $\le n$ and we have for $V_1= \ldots = V_i=R/I$:
\[\rmcr_i(G)(R/I, \ldots, R/I) \cong 
\oplusm_{{(n_1, \ldots, n_i) \in \{1, \ldots, n\}^i \atop n_1 + \ldots + n_i =n}}
B_k\left((R/I \,\, \stackrel{\id}{\ra}\,\, R/I)^{\otimes i}\right) \otimes
(I/I^2)^{\otimes k}.\]
Thus, for all $i\ge 1$, we obtain an isomorphism
\[\alpha_i(R/I, \ldots, R/I): \rmcr_i(F)(R/I, \ldots, R/I) \,\,
\tilde{\ra}\,\, \rmcr_i(G)(R/I, \ldots, R/I).\]
The multiplication with $\bar{r}\in R/I$ in the $l$-th component of
$\rmcr_i(F)(R/I, \ldots, R/I)$ obviously corresponds to the endomorphism
\[\oplusm_{(n_1, \ldots, n_i)} r^{n_l}\quad  {\rm of} \quad
\oplusm_{(n_1, \ldots, n_i)} 
B_k\left((R\ra R)^{\otimes i}\right) \otimes (I/I^2)^{\otimes k}.\] 
The analogous
statement also holds for $G$. Thus, the isomorphisms
$\alpha_i(R/I, \ldots, R/I)$, $i=1, \ldots, n$, are compatible with the action
of $R/I$ in each component. For any $\varepsilon \in \{1,\ldots, n\}^i$ and
$V_1, \ldots, V_i \in \cP_{R/I}$, the diagonal map
\[\gD_{\varepsilon}: \rmcr_i(F)(V_1, \ldots, V_i) \ra 
\rmcr_{|\varepsilon|}(F)(V_1, \ldots, V_1, \: \ldots, \: V_i, \ldots, V_i)\]
is induced by an obvious complex homomorphism
\[\gD_{\varepsilon}: \oplusm_{(n_1, \ldots, n_i)}
\Kos^{n_1}(f_{V_1}) \otimes \ldots \otimes \Kos^{n_i}(f_{V_i})
\ra \oplusm_{(n_1, \ldots, n_{|\varepsilon|})} 
\Kos^{n_1}(f_{V_1}) \otimes \ldots \otimes \Kos^{n_{|\varepsilon|}}(f_{V_i}).\]
Similarly, the diagonal map 
\[\gD_{\varepsilon}: \rmcr_i(G)(V_1, \ldots, V_i) \ra 
\rmcr_{|\varepsilon|}(G)(V_1, \ldots, V_1, \: \ldots, \: V_i, \ldots, V_i)\]
is induced by an obvious complex homomorphism
\[\gD_\varepsilon: 
\oplusm_{(n_1, \ldots, n_i)} \Kos^{n_1}(\id_{V_1}) \otimes \ldots \otimes
\Kos^{n_i}(\id_{V_i}) \ra
\oplusm_{(n_1, \ldots, n_{|\varepsilon|})} \Kos^{n_1}(\id_{V_1}) \otimes \ldots \otimes
\Kos^{n_{|\varepsilon|}}(\id_{V_i}).\]
Using these facts, one easily sees that the isomorphisms $\alpha_i(R/I, \ldots, R/I)$,
$i\ge 1$, are compatible with the diagonal maps. In the same way one can prove that
these isomorphisms are compatible with the plus maps. Hence, by Theorem 1.5, 
the functors $F$ and $G$ are isomorphic, and Theorem 3.2 is proved.

In the next remark, we explicitly describe the isomorphism of Theorem 3.2.

{\bf Remark 3.4.} The isomorphism
\[L_k^n(V) \otimes I^k/I^{k+1} \,\, \tilde{\ra} \,\,
H_k(\Kos^n(f))\]
constructed in the proof of Theorem 3.2 maps the element
\[d_{k+1}(v_1\wedge \ldots \wedge v_{k+1} \otimes v_{k+2} \cdots v_n) \otimes
\overline{r_1\cdots r_k}\] 
to the homology class of the element
\[\sum_{i=1}^{k+1}(-1)^{k+1-i} (r_1q_1) \wedge \ldots \wedge \hat{q}_i \wedge
\ldots \wedge (r_kq_{k+1}) \otimes q_iq_{k+2} \cdots q_n;\] 
here, $q_i$ is a 
preimage of $v_i$ under the given surjective map $Q\ra V$ (for $i=1, \ldots, n$)
and, for any $r\in I$ and $q \in Q$, the element $rq$ of $Q$ is considered as an
element of $P$. 

{\bf Proof.} The isomorphism $G \,\, \tilde{\ra}\,\, F$ constructed in the proof
of Theorem 3.2 induces an isomorphism $\rmcr_n(G) \,\,\tilde{\ra}\,\, 
\rmcr_n(F)$ of $n$-functors. Obviously, the diagrams
\[\xymatrix{
\rmcr_n\left(\gL^{k+1} \otimes \Sym^{n-k-1}\right)(V, \ldots, V) \otimes
I^k/I^{k+1} 
\ar[rrr]^-{\rmcr_n(d_{k+1})(V, \ldots, V) \otimes \id}
\ar[d]^-{+_n \otimes \id} &&&
\rmcr_n(G)(V, \ldots, V)
\ar[d]^-{+_n} \\
\gL^{k+1}(V) \otimes \Sym^{n-k-1}(V) \otimes I^k/I^{k+1}
\ar@{->>}[rrr]^-{d_{k+1} \otimes \id} &&&
G(V) 
}\]
and
\[\xymatrix{
\rmcr_n(G)(V, \ldots, V)
\ar[r]^-{\sim} \ar[d]^-{+_n}& 
\rmcr_n(F)(V, \ldots, V) \ar[d]^-{+_n}\\
G(V) \ar[r]^-{\sim} & F(V)
}\]
commute for all $V\in \cP_{R/I}$. The $R/I$-module 
$\rmcr_n\left(\gL^{k+1}\otimes \Sym^{n-k-1}\right)(V, \ldots, V)$ is
isomorphic to a direct sum of ${n \choose k+1}$ copies of $V^{\otimes n}$. We index
these summands by the subsets $T$ of $\{1, \ldots, n\}$ with $k+1$ elements in the
canonical way. For $v_1 \otimes \ldots \otimes v_n \in V^{\otimes n}$ and
$T=\{t_1 < \ldots < t_{k+1}\}$, let $(v_1 \otimes \ldots \otimes v_n)[T]$ denote
the corresponding element of $\rmcr_n\left(\gL^{k+1} \otimes \Sym^{n-k-1}\right)
(V, \ldots, V)$. We obviously have
\[+_n((v_1 \otimes \ldots \otimes v_n)[T]) = v_{t_1} \wedge \ldots \wedge
v_{t_{k+1}} \otimes \prod_{t\in \{1, \ldots, n\}\backslash T} v_t \quad
{\rm in} \quad \gL^{k+1}(V) \otimes \Sym^{n-k-1}(V).\]
In particular, the restriction of $+_n$ to the direct summand which corresponds
e.g.\ to $\{1, \ldots, k+1\}$ is surjective. Thus it suffices to show that
the composition of the upper horizontal maps of the above diagrams 
with the right vertical map in
the lower diagram maps the element $(v_1 \otimes \ldots \otimes v_n)
[\{1, \ldots, k+1\}] \otimes \overline{r_1\cdots r_k}$ to the homology class
of the element given in Remark 3.4. We obviously have
\begin{eqnarray*}
\lefteqn{\rmcr_n(d_{k+1})(V, \ldots, V)((v_1 \otimes \ldots \otimes v_n)
[\{1, \ldots, k+1\}])}\\
&& = \sum_{i=1}^{k+1} (-1)^{k+1-i} 
(v_1 \otimes \ldots \otimes v_n)[\{1, \ldots, \hat{i}, \ldots, k+1\}]
\end{eqnarray*}
in $\rmcr_n\left(\gL^k \otimes \Sym^{n-k}\right)(V, \ldots, V)$. The
isomorphism 
\[\rmcr_n(G)(V, \ldots, V) \,\, \tilde{\ra}\,\, \rmcr_n(F)(V, \ldots, V)\]
is equal to the composition
\begin{eqnarray*}
\lefteqn{\rmcr_n(G)(V, \ldots, V) \cong V^{\otimes n} \otimes \rmcr_n(G)(R/I, \ldots, R/I)}\\
&&\cong V^{\otimes n} \otimes \rmcr_n(F)(R/I, \ldots, R/I) \cong 
\rmcr_n(F)(V, \ldots, V)
\end{eqnarray*}
where the middle isomorphism has been explicitly constructed in the proof
of Theorem 3.2 and the exterior isomorphisms have been introduced in Example 1.7.
The composition 
\begin{eqnarray*}
\lefteqn{\Image\left(\rmcr_n(d_{k+1})(R/I, \ldots, R/I)\right) \otimes I^k/I^{k+1} }\\
&& \cong \rmcr_n(G)(R/I, \ldots, R/I) \cong \rmcr_n(F)(R/I, \ldots, R/I)
\end{eqnarray*}
maps the element $\rmcr_n(d_{k+1})(R/I, \ldots, R/I)(1[\{1, \ldots, k+1\}])
\otimes \overline{r_1\cdots r_k}$ to the homology class of the element
\[\sum_{i=1}^{k+1} (-1)^{k+1-i}(r_1 \otimes \ldots \otimes r_{i-1} \otimes 1 \otimes
r_i \otimes \ldots \otimes r_k \otimes 1 \otimes \ldots \otimes 1)
[\{1, \ldots, \hat{i}, \ldots, k+1\}]\]
in the complex $(I\ra R)^{\otimes n}$
(see the proof of Theorem 3.2). Furthermore, the isomorphism 
$V^{\otimes n} \otimes \rmcr_n(F)(R/I, \ldots, R/I) \cong 
\rmcr_n(F)(V, \ldots, V)$ maps the element $v_1 \otimes \ldots \otimes v_n \otimes x$
to the image of $x$ under the map
\[\rmcr_n(F)(v_1, \ldots, v_n): \rmcr_n(F)(R/I, \ldots, R/I) \ra 
\rmcr_n(F)(V, \ldots, V).\]
(This easily follows from Example 1.7). Hence, the composition of the
upper horizontal maps in the above diagrams maps the element
$(v_1 \otimes \ldots \otimes v_n)[\{1, \ldots, k+1\}]\otimes \overline{r_1 \cdots r_k}$
to the homology class of the element
\begin{eqnarray*}
\lefteqn{\sum_{i=1}^{k+1}(-1)^{k+1-i}((r_1q_1) \otimes \ldots \otimes
(r_{i-1}q_{i-1}) \otimes q_i }\\
&& \otimes (r_iq_{i+1}) \otimes \ldots \otimes
(r_kq_{k+1}) \otimes q_{k+2} \otimes \ldots \otimes q_n)[\{1, \ldots, \hat{i}, 
\ldots, k+1\}]
\end{eqnarray*}
in the complex $(P\ra Q)^{\otimes n}$. The plus map $+_n: 
\rmcr_n(F)(V, \ldots, V) \ra F(V)$ maps the latter homology class to the homology class
stated in Remark 3.4. Thus, Remark 3.4 is proved.

In the following remark, we sketch an alternative proof of Theorem 3.2.

{\bf Remark 3.5.} An easy, but tedious calculation shows directly that the association 
described in Remark 3.4 yields a well-defined homomorphism
\[G(V) = L_k^n(V) \otimes I^k/I^{k+1} \ra H_k(\Kos^n(f)) = F(V)\]
for all $V\in \cP_{R/I}$. Furthermore, one can verify that the induced
homomorphisms $\rmcr_i(G)(R/I, \ldots, R/I) \ra \rmcr_i(F)(R/I, \ldots, R/I)$, 
$i\ge 1$, are bijective. By Proposition 1.2, this implies that the morphism
$G \ra F$ of functors is an isomorphism.\\
So, in this proof, the abstract argument Theorem 1.5 has been replaced by more explicit
computations which, from the combinatorial point of view, are more complicated.

{\bf Remark 3.6.} Analogously to Theorem 3.2, one can prove the following statement
for the Koszul complex $\tilde{\Kos}^n(f)$ introduced in Remark 2.7:
Under the assumptions of Theorem 3.2, we have:
\[H_k(\tilde{\Kos}^n(f)) \cong \tilde{L}^n_k(V) \otimes (I/I^2)^{\otimes k}
\quad \textrm{for all }k\ge 0;\]
here, $\tilde{L}^n_k(V)$ denotes the coSchur module of type $(k+1, 1, \ldots, 1)$
which can be defined analogously to 3.1 (see also Definition II 1.3 on p.\ 220
in \cite{ABW}). 

\bigskip

\section*{\S 4 A New Proof of the Adams-Riemann-Roch Formula for Regular Closed
Immersions}

In this section, we deduce the Adams-Riemann-Roch theorem for regular closed
immersions of codimension $1$ from the main result of the previous section.
Using blowing up and the excess intersection formula, the Adams-Riemann-Roch
theorem for regular closed immersions of arbitrary codimension then follows 
from this as in \cite{SGA6}. In particular, we obtain an easy, elementary
and natural proof of the Adams-Riemann-Roch theorem which does not use
the deformation to the normal cone any longer.

First, we recall the following definitions and facts of the theory of
$\lambda$-rings (see Expose V in \cite{SGA6}, section 4.5 in \cite{So}, or
Chapter I in \cite{FL}). Let $K$ be a (special) $\lambda$-ring, $C \in K$
of (finite) $\lambda$-degree $d$, and $x \in K$ arbitrary. Furthermore, let
$\mu$ be a polynomial in the $\lambda$-operations $\lambda_1, \lambda_2, \ldots$
without constant term. Then, there is a unique element $\mu(C,x) \in K$ which
is functorial in $(K,C,x)$ and which has the following property:
\[\mu(x\cdot \lambda_{-1}(C)) = \mu(C,x) \cdot \lambda_{-1}(C)\]
(here, $\lambda_{-1}(C) := \sum_{i=0}^d (-1)^i \lambda_i(C)$). This immediately
follows from the following fact which is easy to prove (cf.\ Lemma 5.2 in
Expose V of \cite{SGA6}): If the elements $\lambda_1(C), \ldots, \lambda_d(C)$
and $\lambda_1(x), \lambda_2(x), \ldots$ are considered as indeterminates
and if the polynomial ring $\ZZ[\lambda_1(C), \ldots, \lambda_d(C), 
\lambda_1(x), \lambda_2(x), \ldots]$ is equipped with the obvious 
$\lambda$-structure, then the element $\mu(x \cdot \lambda_{-1}(C))$ is
divisible by $\lambda_{-1}(C)$. For $n\ge 1$, let $\psi_n$ denote the 
$n$-th {\em Adams operation}, i.e., $\psi_n = N_n(\lambda_1, \ldots, \lambda_n)$
where $N_n$ denotes the $n$-th Newton polynomial. Furthermore, let 
$\theta^n(C)$ denote the $n$-th Bott element of $C$ (see \S 4 in
\cite{KoGRR} or p.\ 24 in \cite{FL}). The $n$-th {\em symmetric power operation}
is inductively defined by $\gs_0 \equiv 1$ and $\gs_n = \sum_{i=1}^{n}
(-1)^{i-1} \lambda_i \gs_{n-i}$ for $n\ge 1$. For $0\le k \le n$, the
{\em Schur operation} $s_k^n$ {\em of type} $(k+1, 1, \ldots, 1)$ is defined by
$s_k^n = \sum_{i=k+1}^n (-1)^{i-k-1} \lambda_i \gs_{n-i}$.

{\bf Lemma 4.1.} \\
(a) Let $\mu$ be the product of two polynomials $\mu_1$ and $\mu_2$. Then we have:
\[\mu(C,x) = \mu_1(C,x) \cdot \mu_2(C,x) \cdot \lambda_{-1}(C).\]
(b) For all $n\ge 0$ and $x,y \in K$, we have:
\[\gs_n(C, x+y) = \gs_n(C,x) + \sum_{i=1}^{n-1} \gs_i(C,x) \cdot
\gs_{n-i}(C,y) \cdot \lambda_{-1}(C) + \gs_n(C,y).\]
(c) For all $n\ge 1$ we have: $\psi_n(C,x) = \theta^n(C) \cdot \psi_n(x)$.\\
(d) Let $d=1$. Then we have for all $n\ge 1$:
\[\gs_n(C,x) = \sum_{k=0}^{n-1} (-1)^k s_k^n(x) C^k.\]
In particular, we have $\gs_n(1,x) = \psi_n(x)$.

{\bf Proof.} The assertions (a) and (b) immediately follow from the definition.
The assertion (c) follows from the multiplicativity of $\psi_n$. If $d=1$, then
we have:
\begin{eqnarray*}
\lefteqn{\gs_n(x \cdot \lambda_{-1}(C)) = 
\gs_n(x \cdot (1-C)) = \sum_{i=0}^n \gs_{n-i}(x) \cdot \gs_i(-xC) }\\
&& = \sum_{i=0}^n (-1)^i \gs_{n-i}(x) \cdot \lambda_i(xC) \qquad (\textrm{since }
\gs_i(-x) = (-1)^i\lambda_i(x))\\
&& =\sum_{i=0}^n (-1)^i \gs_{n-i}(x) \cdot \lambda_i(x) \cdot C^i \qquad
(\textrm{since } C \textrm{ is of } \lambda\textrm{-degree } 1)\\
&&= \sum_{i=0}^n (-1)^i \gs_{n-i}(x) \cdot \lambda_i(x) \cdot (C^i -1) \qquad
(\textrm{since } \sum_{i=0}^{n} (-1)^i \gs_{n-i}(x) \cdot \lambda_i(x) =0)\\
&& = \sum_{i=0}^n (-1)^{i+1} \gs_{n-i}(x) \cdot \lambda_i(x)\cdot (1 + C + \ldots + C^{i-1})
\cdot \lambda_{-1}(C) \quad (\textrm{geometric sum}).
\end{eqnarray*}
Hence we have:
\begin{eqnarray*}
\lefteqn{\gs_n(C,x) = \sum_{i=1}^n (-1)^{i+1} \gs_{n-i}(x)\cdot \lambda_i(x)\cdot
(1+C + \ldots + C^{i-1})}\\
&& = \sum_{k=0}^{n-1}\left(\sum_{i=k+1}^n (-1)^{i+1} 
\gs_{n-i}(x) \cdot \lambda_i(x)\right)C^k
= \sum_{k=0}^{n-1} (-1)^k s_k^n(x) C^k.
\end{eqnarray*}
In particular, we have $\gs_n(1,x) = \sum_{k=0}^{n-1} (-1)^k s_k^n(x) = \psi_n(x)$
by section 3 in \cite{Gr}. 

Now, let $X$ be a noetherian scheme with the property that each coherent $\cO_X$-module is
a quotient of a locally free $\cO_X$-module (of finite type). Let $K_0(X)$ denote
the Grothendieck group of all locally free $\cO_X$-modules and $K_0^\infty(X)$ 
the Grothendieck group of all coherent $\cO_X$-modules which have a finite resolution
by locally free $\cO_X$-modules. By Proposition 4.1 on p.\ 126 in \cite{FL}, the
canonical homomorphism $K_0(X) \ra K_0^\infty(X)$ is bijective. Let $i:Y \hookrightarrow
X$ be a regular closed immersion of codimension $d$ with conormal sheaf $\cC$. 
Then, for any locally free $\cO_Y$-module $\cV$, the direct image $i_*(\cV)$ is 
a coherent $\cO_X$-module which has a finite locally free resolution (see 
p.\ 127 in \cite{FL}). Thus, the association $[\cV] \mapsto [i_*(\cV)]$ induced
a well-defined homomorphism
\[i_*: K_0(Y) \ra K_0^\infty(X) \cong K_0(X).\]
The following formula describes the behavior of products with respect to
this homomorphism. It follows from the self intersection formula and it
will be used in our proof of the Adams-Riemann-Roch formula.

{\bf Proposition 4.2.} For all $x,y \in K_0(Y)$, we have:
\[i_*(x) \cdot i_*(y) = i_*(x \cdot y \cdot \lambda_{-1}([\cC])) \quad
\textrm{in} \quad K_0(X).\]

{\bf Proof.} See Corollaire 2.8 on p.\ 436 in \cite{SGA6}.

The assertion (a) of the following theorem has already been proved in
Th\'eor\`eme 4.3 on p.\ 449 in \cite{SGA6}, however only modulo torsion. The
claimed version without denominators has been proved in Th\'eor\`eme 2.1
on p.\ 24 in \cite{J} using the deformation to the normal cone.  The 
Adams-Riemann-Roch formula given in (b) follows also from
Th\'eor\`eme 2.1 on p.\ 24 in \cite{J} and it has been proved in
Theorem 6.3 on p.\ 142 in \cite{FL}.

{\bf Theorem 4.3.} Let $n\ge 1$.\\
(a) (Riemann-Roch formula for $\gs_n$) For all $y \in K_0(Y)$ we have:
\[\gs_n(i_*(y))=i_*(\gs_n([\cC],y)) \quad \textrm{in} \quad K_0(X).\]
(b) (Adams-Riemann-Roch formula) For all $y \in K_0(Y)$ we have:
\[\psi_n(i_*(y)) = i_*(\theta^n([\cC]) \cdot \psi_n(y)) \quad \textrm{in}
\quad K_0(X).\]

{\bf Proof.} \\
(a) First, let $d=1$. By Lemma 4.1(b) and Proposition 4.2, the elements $y$ of
$K_0(Y)$ for which the Riemann-Roch formula for $\gs_n$ holds for all $n\ge 1$
is a subgroup of $K_0(Y)$. Thus, it suffices to prove the Riemann-Roch 
formula for $y=[\cV]$ where $\cV$ is a locally free $\cO_Y$-module. We choose
a locally free resolution $0\ra \cP \,\, \stackrel{f}{\ra} \,\,\cQ 
\ra i_*(\cV) \ra 0$ of $i_*(\cV)$ on
$X$. By gluing the isomorphisms of Theorem 3.2, we obtain isomorphisms
\[\cH_k(\Kos^n(f)) \cong i_*(L_k^n(\cV)\otimes \cC^{\otimes k}), \quad 
k=0, \ldots, n.\]
Hence, we have in $K_0^\infty(X)$:
\[\sum_{k=0}^n (-1)^k[\gL^k(\cP) \otimes \Sym^{n-k}(\cQ)] = 
i_*\left(\sum_{k=0}^{n-1} (-1)^k [L_k^n(\cV) \otimes \cC^{\otimes k}]\right).\]
Thus, we obtain:
\begin{eqnarray*}
\lefteqn{\gs_n(i_*([\cV])) = \gs_n([\cQ]-[\cP]) = 
\sum_{k=0}^{k} (-1)^k [\gL^k(\cP)\otimes \Sym^{n-k}(\cQ)]}\\
&&=i_*\left(\sum_{k=0}^{n-1} (-1)^k s_k^n([\cV]) \cdot [\cC]^k\right)
= i_*(\gs_n([\cC],[\cV])) \quad \textrm{in} \quad K_0(X);
\end{eqnarray*}
here, we have used the exactness of the Koszul complex $\Kos^n(\id_\cV)$ and
Lemma 4.1(d). This proves assertion (a) in the case $d=1$. The general case
can be deduced from this as on p.\ 449 and p.\ 450 in \cite{SGA6} using blowing up and
the excess intersection formula.\\
(b) Let $\psi_n = \sum_\nu a_\nu \gs_1^{\nu_1} \cdots \gs_n^{\nu_n}$ be the 
representation of $\psi_n$ as a polynomial in the symmetric power operations
$\gs_1, \ldots, \gs_n$. Then we have in $K_0(X)$:
\begin{eqnarray*}
\lefteqn{\psi_n(i_*(y)) = \sum_\nu a_\nu \gs_1(i_*(y))^{\nu_1} \cdots
\gs_n(i_*(y))^{\nu_n}}\\
&& = \sum_\nu a_\nu \, i_*(\gs_1([\cC],y))^{\nu_1} \cdots 
i_*(\gs_n([\cC],y))^{\nu_n} \quad (\textrm{Theorem 4.3(a)})\\
&& =i_*\left(\sum_\nu a_\nu \, \gs_1([\cC],y)^{\nu_1} \cdots
\gs_n([\cC],y)^{\nu_n} \cdot \lambda_{-1}([\cC])^{\nu_1 + \ldots + \nu_n-1}\right)
\quad (\textrm{Proposition 4.2})\\
&& = i_*(\psi_n([\cC],y)) \qquad \qquad (\textrm{Lemma 4.1(a)})\\
&& = i_*(\theta^n([\cC]) \cdot \psi_n(y)) \qquad \qquad (\textrm{Lemma 4.1(c)}).
\end{eqnarray*}
This ends the proof of Theorem 4.3.

\bigskip

\section*{\S 5 Riemann-Roch for Tensor Powers}

In the paper \cite{KoTe}, we have established Riemann-Roch formulas for tensor
power operations. In the case of regular closed immersions, we have used the 
deformation to the normal cone to prove them. The aim of this section is to
prove these formulas by using the method developed in the previous sections.
In contrast to the previous sections, this method can here be applied to
regular closed immersions not only of codimension $1$ but of arbitrary
codimension.

Let $R$ be a noetherian ring and $n\in \NN$. Let $\gS_n$ denote the
$n$-th symmetric group, i.e., the group of permutations of $I_n:=\{1, \ldots, n\}$. 
We write $R[I_n]$ for the direct sum of $n$ copies of $R$ together with the 
obvious action of $\gS_n$. The $\gS_n$-module $K_n$ is defined by the 
short exact sequence
\[0\longrightarrow K_n \longrightarrow R[I_n] \,\, 
\stackrel{{\rm sum}}{\longrightarrow} R \longrightarrow 0.\]
Let $I$ be an ideal in $R$ which is locally generated by a regular sequence of 
length $d$. Furthermore, let $P. \ra V$ be an $R$-projective resolution of
a f.\ g.\ projective $R/I$-module $V$. We view the total complex $P.^{\otimes n}$ of the 
$n$-th tensor power of the complex $P.$ as a complex of $\gS_n$-modules in the
canonical way. In particular, the homology modules $H_k(P.^{\otimes n})$, $k\ge 0$,
are equipped with a canonical $\gS_n$-action.

{\bf Theorem 5.1.} For all $k\ge 0$, there is a canonical $\gS_n$-isomorphism
\[H_k(P.^{\otimes n}) \cong V^{\otimes n} \otimes \gL^k(I/I^2 \otimes K_n);\]
here, $V^{\otimes n}$ is considered as a $\gS_n$-module with the natural action
and $I/I^2$ as a $\gS_n$-module 
with the trivial action. In particular, $H_k(P.^{\otimes n})$ is a
projective $R/I$-module.

{\bf Proof.} We fix $k\in \{0, \ldots, n\}$. A homotopy between two complex
homomorphisms $\xymatrix@1{f,g: P. \ar@<0.5ex>[r] \ar@<-0.5ex>[r] & Q.}$ induces
a (non-equivariant) homotopy between the complex homomorphisms
$\xymatrix@1{f^{\otimes n}, g^{\otimes n}: P.^{\otimes n} \ar@<0.5ex>[r] \ar@<-0.5ex>[r]
& Q.^{\otimes n}}$. In particular, the $\gS_n$-homomorphisms $H_k(f^{\otimes n})$
and $H_k(g^{\otimes n})$ coincide. Therefore, we obtain a well-defined
functor
\[F: \cP_{R/I} \ra (R[\gS_n]\textrm{-modules}), \quad V \mapsto H_k(P.^{\otimes n}),\]
as in the proof of Theorem 3.2. Let $G$ denote the functor
\[G: \cP_{R/I} \ra (R[\gS_n]\textrm{-modules}), \quad V \mapsto 
V^{\otimes n} \otimes \gL^k(I/I^2 \otimes K_n).\]
We have to show that $F(V) \cong G(V)$ for all $V \in \cP_{R/I}$. Again, we verify
the assumptions of Theorem 1.5.\\
First, we suppose that the ideal $I$ is globally generated by a regular sequence.
This means that there is a homomorphism $\varepsilon: F\ra R$ from a free
$R$-module $F$ of rank $d$ to $R$ such that the Koszul complex $\Kos^d(\varepsilon)$ is
a resolution of $R/I$ (by virtue of the canonical map $\Kos^d(\varepsilon)_0 \cong
R \rightepi R/I$). One easily checks that the complex $(\Kos^d(\varepsilon))^{\otimes n}$
is $\gS_n$-isomorphic to the Koszul complex $\Kos^d(F\otimes R[I_n] 
\,\,\stackrel{{\rm sum}}{\longrightarrow}\,\, F 
\,\, \stackrel{{\varepsilon}}{\longrightarrow}\,\, R)$. By Lemma 3.3 in \cite{KoTe},
the homology module $H_k(\Kos^d(F\otimes R[I_n] \ra F \ra R))$ is $\gS_n$-isomorphic
to the $k$-th exterior power of 
\[\textrm{ker}(F/IF \otimes R[I_n] 
\,\, \stackrel{{\rm sum}}{\longrightarrow} \,\, F/IF
\,\,\stackrel{\bar{\varepsilon}}{\longrightarrow}\,\, I/I^2) \cong 
F/IF \otimes K_n \cong I/I^2 \otimes K_n.\]
Thus, we have established a $\gS_n$-isomorphism $F(R/I) \cong \gL^k(I/I^2 \otimes K_n)$.
It does obviously not depend on the chosen homomorphism $\varepsilon$. 
Hence, we can glue these isomorphisms and obtain a canonical 
$\gS_n$-isomorphism $F(R/I) \cong \gL^k(I/I^2 \otimes K_n)$
also in the case the ideal $I$ is only locally generated by
a regular sequence. For any two complexes $P.$ and $Q.$, we have:
\[(P. \oplus Q.)^{\otimes n} \cong \oplusm _{n_1+n_2=n} \Ind_{\gS_{n_1} \times \gS_{n_2}}
^{\gS_n} (P.^{\otimes n_1} \otimes Q.^{\otimes n_2}).\]
Hence, as in the proof of Theorem 3.2, we obtain for all $i\ge 0$:
\begin{eqnarray*}
\rmcr_i(F)(R/I, \ldots, R/I) &\cong &
\oplusm_{{(n_1, \ldots, n_i) \in \{1, \ldots, n\}^i \atop n_1 + \ldots +n_i =n}}
\Ind_{\gS_{n_1} \times \ldots \times \gS_{n_i}}^{\gS_n}
\Res_{\gS_{n_1} \times \ldots \times \gS_{n_i}}^{\gS_n} F(R/I)\\
& \cong & \left(\oplusm_{(n_1, \ldots, n_i)} 
\Ind_{\gS_{n_1} \times \ldots \times \gS_{n_i}}^{\gS_n}(R) \right)
\otimes \gL^k(I/I^2 \otimes K_n) \\
& \cong & \rmcr_i(G)(R/I, \ldots, R/I).
\end{eqnarray*}
One easily sees that these $\gS_n$-isomorphisms are compatible with the action
of $R/I$ in each component and with the diagonal and plus maps. Thus, Theorem
5.1 follows from Theorem 1.5.

{\bf Example 5.2.} Let $n=2$. Then Theorem 5.1 yields the well known 
$\gS_2$-isomorphism
\[{\rm Tor}_k^R(V,V) \cong V^{\otimes 2} \otimes \gL^k((I/I^2)_\sgn)\]
where the transposition in $\gS_2$ acts by $-1$ on $(I/I^2)_\sgn$.

Now, let $X$ be a noetherian scheme with the property that each coherent $\cO_X$-module is
a quotient of a locally free $\cO_X$-module. Let $K_0(\gS_n,X)$ denote the 
Grothendieck group of all locally free $\gS_n$-modules on $X$. Using the
``binomial theorem''
\[(\cF \oplus \cG)^{\otimes n} \cong \oplusm_{n_1 + n_2 = n}
\Ind_{\gS_{n_1} \times \gS_{n_2}}^{\gS_n} (\cF^{\otimes n_1} \otimes
\cG^{\otimes n_2}),\]
the association $[\cF] \mapsto [\cF^{\otimes n}]$ can be canonically extended 
to a map
\[\tau_n: K_0(X) \ra K_0(\gS_n,X)\]
(see section 1 in \cite{KoTe}). It is called the $n$-th tensor power operation.
Let $i:Y\hookrightarrow X$ be a regular closed immersion of codimension $d$
with conormal sheaf $\cC$. Let $i_*$ denote both the homomorphism from
$K_0(Y)$ to $K_0(X)$ introduced in section 4 and the homomorphism from
$K_0(\gS_n,Y)$ to $K_0(\gS_n,X)$ defined analogously. The following theorem is
Theorem 4.1 of \cite{KoTe}.

{\bf Theorem 5.3} (Riemann-Roch formula for $\tau_n$). For all $y \in K_0(Y)$,
we have:
\[\tau_n(i_*(y)) = i_*(\lambda_{-1}([\cC \otimes K_n]) \cdot \tau_n(y)) \quad
\textrm{in} \quad K_0(\gS_n,X).\]

{\bf Proof.} If Theorem 5.3 is true for $y_1, y_2 \in K_0(Y)$ and for all $n\ge 1$,
then it is also true for $y_1+y_2$. This follows from the following calculation:
\begin{eqnarray*}
\lefteqn{\tau_n(i_*(y_1 +y_2)) = \sum_{n_1+n_2 =n}
\Ind_{\gS_{n_1}\times \gS_{n_2}}^{\gS_n}
\Big(\tau_{n_1}(i_*(y_1)) \cdot \tau_{n_2}(i_*(y_2))\Big)}\\
&& =  \sum_{n_1+n_2 =n} \Ind_{\gS_{n_1}\times \gS_{n_2}}^{\gS_n}
\Big(i_*(\lambda_{-1}([\cC \otimes K_{n_1}]) \cdot \tau_{n_1}(y_1)) \cdot
i_*(\lambda_{-1}([\cC \otimes K_{n_2}]) \cdot \tau_{n_2}(y_2)) \Big)\\
&& =  \sum_{n_1+n_2 =n} \Ind_{\gS_{n_1}\times \gS_{n_2}}^{\gS_n}
\Big(i_*\Big(\lambda_{-1}([\cC \otimes K_{n_1}])  
\lambda_{-1}([\cC \otimes K_{n_2}]) \lambda_{-1}([\cC]) \tau_{n_1}(y_1)
\tau_{n_2}(y_2) \Big) \Big)\\
&& = \sum_{n_1+n_2 =n} \Ind_{\gS_{n_1}\times \gS_{n_2}}^{\gS_n}
\Big(i_*\Big(\Res_{\gS_{n_1} \times \gS_{n_2}}^{\gS_n}
\lambda_{-1}([\cC \otimes K_n]) \cdot \tau_{n_1}(y_1) \cdot \tau_{n_2}(y_2) \Big)\Big)\\
&& = i_*\left(\lambda_{-1}([\cC \otimes K_n]) \cdot
\sum_{n_1 + n_2=n} \Ind_{\gS_{n_1} \times \gS_{n_2}}^{\gS_n}
\tau_{n_1}(y_1) \cdot \tau_{n_2}(y_2)\right)\\
&& = i_*\left(\lambda_{-1}([\cC \otimes K_n]) \cdot \tau_n(y_1+y_2)\right).
\end{eqnarray*}
(Here, we have used Proposition 4.2 and Frobenius reciprocity.) Thus, it suffices
to prove Theorem 5.3 for $y = [\cV]$ where $\cV$ is a locally free $\cO_Y$-module.
Let $\cP. \ra i_*(\cV)$ be a locally free resolution of $i_*(\cV)$ on $X$. By
gluing the isomorphisms constructed in Theorem 5.1, we obtain $\gS_n$-isomorphisms
\[\cH_k(\cP.^{\otimes n}) \cong i_*(\cV^{\otimes n} \otimes
\gL^k(\cC \otimes K_n)), \quad k\ge 0.\]
Hence, we have in $K_0(\gS_n,X)$:
\begin{eqnarray*}
\lefteqn{\tau_n(i_*([\cV])) = \tau_n\left(\sum_{i\ge 0} (-1)^k[\cP_k]\right)
= \sum_{k\ge 0} (-1)^k \left[(\cP.^{\otimes n})_k\right]}\\
&&= \sum_{k\ge 0} (-1)^k i_*\left([\cV^{\otimes n}\otimes \gL^k(\cC \otimes K_n)]\right)
= i_*\left(\lambda_{-1}([\cC\otimes K_n]) \cdot \tau_n([\cV])\right).\\
\end{eqnarray*}
This proves Theorem 5.3.

\bigskip

\section*{\S 6 Some Computations for Regular Immersions of Higher Codimension}

Let $X$ be a noetherian scheme with the property 
that each coherent $\cO_X$-module is a quotient
of a locally free $\cO_X$-module. Let $i: Y \hookrightarrow X$ be a regular closed immersion
of codimension $d$ with conormal sheaf $\cC$ and $n\in \NN$. In this section,
we study the question whether, in the case $d\ge 2$, 
the Riemann-Roch formula for $\gs_n$ and $i_*$ (cf.\ Theorem 4.3(a)) can be 
proved directly with the methods developed in the previous sections (i.e., 
without using blowing up). More precisely, for any locally free $\cO_Y$-module
$\cV$ and any locally free resolution $\cP. \ra i_*(\cV)$ of $i_*(\cV)$ on $X$,
we are looking for a complex $\cK.^n$ of locally free $\cO_X$-modules such that
we have 
\[\gs_n(i_*([\cV])) \,\, \stackrel{{\rm def}}{=} \,\, 
\gs_n(\sum_{k\ge 0} (-1)^k [\cP_k]) = \sum_{k\ge 0}(-1)^k [\cK_k^n] \quad
\textrm{in} \quad K_0(X)\] 
and 
\[\gs_n([\cC], [\cV]) = \sum_{k\ge 0} (-1)^k [\cH_k(\cK.^n)] \quad \textrm{in}
\quad K_0(Y)\] 
(or at least in the relative $K_0$-group of $Y$ in $X$, see \cite{TT}).
The existence of such a complex would immediately imply the Riemann-Roch formula
for $\gs_n$ (see the proof of Theorem 4.3(a)). In the case $d=1$, we have
used the complex $\Kos^n(\cP_1 \ra \cP_0)$ (see section 4) which is 
quasi-isomorphic to the complex $N\Sym^n\gG(\cP_1 \ra \cP_0)$ by Proposition 2.4. In 
section 5, we have used the complex $\cP.^{\otimes n}$ which is quasi-isomorphic
to the complex $N(\gG(\cP.)^{\otimes n})$ by the Eilenberg-Zilber theorem
(see Application 8.5.3 on p.\ 277 in \cite{We}). Thus, a natural candidate 
for the above complex $\cK.^n$ is the complex $\cK.^n := N \Sym^n \gG(\cP.)$.
We prove in Lemma 6.1 that this complex satisfies at least the first
condition $\gs_n(i_*([\cV])) = \sum_{k\ge 0} (-1)^k [\cK_k^n]$. In Lemma 6.2,
we give a representation of $\gs_2([\cC], [\cV])$ as an alternating sum
which suggests the following conjecture: $\cH_k(\cK.^2) \cong 
\Sym^2(\cV) \otimes \gL^k(\cC)$ for $k$ even and $\cH_k(\cK.^2) \cong
\gL^2(\cV) \otimes \gL^k(\cC)$ for $k$ odd. In fact, this conjecture is true
if $2$ is invertible on $X$ (see Corollary 6.3). Theorem 6.4 states that,
in the case $d=2$, we have $\cH_0(\cK.^2) \cong \Sym^2(\cV)$, 
$\cH_1(\cK.^2) \cong \gL^2(\cV) \otimes \cC$, however $\cH_2(\cK.^2) \cong
D^2(\cV) \otimes \gL^2(\cV)$, i.e., the conjecture is false for $k=2$. 
Nevertheless, this theorem implies the second condition for the
complex $\cK.^2$ in the case $d=2$ since $[D^2(\cV)] = [\Sym^2(\cV)]$ in
$K_0(Y)$. In the case $d > 2$, $\cH_k(\cK.^2)$ probably has non-trivial
$2$-torsion for $k\ge 2$ (see Remark 6.5). In particular, it is not clear
how to compute $\cH_k(\cK.^2)$ in general and to check the second condition
for the complex $\cK.^2$. On the other hand, it might be possible to compute 
$\cH_k(\cK.^n)$ in the case $d=2$ even for all $n\ge 2$ and to verify the
second condition for $\cK.^n$ in this case (see Example 6.6). I hope to
say more on this in a future paper.

For any finite complex $\cP.$ of locally free $\cO_X$-modules let
$[\cP.]$ denote the Euler characteristic $\sum_{i\ge 0} (-1)^i [\cP_i] \in K_0(X)$.

{\bf Lemma 6.1.} Let $\cP.$ be a finite complex of locally free $\cO_X$-modules.
Then we have for all $n\ge 0$:
\[[N\Sym^n\gG(\cP.)] = \gs_n([\cP.]) \quad \textrm{in} \quad K_0(X).\]

{\bf Proof.} One easily checks that the complex $N\Sym^n\gG(\cP.)$ has finite
length. Furthermore, the $\cO_X$-modules $N\Sym^n\gG(\cP.)_k$, $k\ge 0$, are
locally free (see section 2 in \cite{JM}). Therefore, the expression
$[N\Sym^n\gG(\cP.)] \in K_0(X)$ is well-defined. First, we suppose that
$\cP.$ is the complex $\cP[-l]$ where $\cP$ is a locally free $\cO_X$-module
and $l\in \NN_0$. We have to show that $\gs_n((-1)^l[\cP]) = 
[N\Sym^n\gG(\cP[-l])]$ in $K_0(X)$. We prove this by induction on $l$. For
$l=0,1$ this follows for example from Proposition 2.4. We write $C(\cP[-l])$ for
the mapping cone of the identity of the complex $\cP[-l]$. The obvious short
exact sequence
\[0 \ra \cP[-l] \ra C(\cP[-l]) \ra \cP[-l-1] \ra 0\]
of complexes induces a natural filtration of the complex $N\Sym^n\gG(C(\cP[-l]))$
whose successive quotients are the complexes $N(\Sym^i\gG(\cP[-l]) \otimes
\Sym^{n-i}\gG(\cP[-l-1]))$, $i=0, \ldots, n$. The complex $N\Sym^n\gG(C(\cP[-l]))$
is homotopic to the zero complex (cf.\ section 2). For all $i=0, \ldots, n$, the
complex $N(\Sym^i\gG(\cP[-l]) \otimes \Sym^{n-i}\gG(\cP[-l-1]))$ is 
quasi-isomorphic to the complex $N\Sym^i\gG(\cP[-l]) \otimes N\Sym^{n-i}\gG
(\cP[-l-1])$ by the Eilenberg-Zilber theorem (see Application 8.5.3 on p.\ 277 in
\cite{We}). Hence, we have for all $n\ge 1$:
\begin{eqnarray*}
&& \sum_{i=0}^n[N\Sym^i\gG(\cP[-l])] \cdot [N\Sym^{n-i}\gG(\cP[-l-1])] \\
& =& \sum_{i=0}^n [N\Sym^i\gG(\cP[-l]) \otimes N\Sym^{n-i}\gG(\cP[-l-1])]\\
&=& \sum_{i=0}^n[N(\Sym^i\gG(\cP[-l]) \otimes \Sym^{n-i}\gG(\cP[-l-1]))]\\
&=& N\Sym^n\gG(C(\cP[-l])) = 0;
\end{eqnarray*}
i.e., the power series $\sum_{i\ge 0} [N\Sym^i\gG(\cP[-l-1])] t^i$ is the inverse of the
power series $\sum_{i \ge 0} [N\Sym^i\gG(\cP[-l])] t^i$. Using the induction hypothesis,
we obtain:
\[[N\Sym^n\gG(\cP[-l-1])] = \gs_n(-[\cP[-l]]) = \gs_n((-1)^{l+1} [\cP])\]
as was to be shown. (Alternatively, this can be deduced from Proposition 4.3.2.1(i) on
p.\ 111 in \cite{Il}.) Now, let $\cP.$ be an arbitrary complex of length $l$. 
We prove Lemma 6.1 by induction on $l$. Let $\gs_{\le l-1}(\cP.)$ denote the
brutal truncation of $\cP.$. The obvious short exact sequence
\[0\ra \gs_{\le l-1}(\cP.) \ra \cP. \ra \cP_l[-l] \ra 0\]
of complexes yields the equality
\[[N\Sym^n\gG(\cP.)] = \sum_{i=0}^n [N\Sym^i\gG(\gs_{\le l-1}(\cP.))] \cdot
[N\Sym^{n-i}\gG(\cP_l[-l])] \quad \textrm{in} \quad K_0(X)\]
as above. Using the induction hypothesis and the formula already proved, we
finally obtain
\[[N\Sym^n\gG(\cP.)] = \sum_{i=0}^n \gs_i \left(\sum_{i=0}^{l-1}(-1)^i [\cP_i]\right)
\cdot \gs_{n-i}((-1)^l[\cP_l]) = \gs_n([\cP.])\]
as was to be shown.

{\bf Lemma 6.2.} Let $K$ be a $\lambda$-ring, $C\in K$ of finite $\lambda$-degree
$d$ and $x \in K$ arbitrary. Then we have:
\[\gs_2(C,x) = \gs_2(x)\lambda_0(C) - \lambda_2(x) \lambda_1(C) +
\gs_2(x) \lambda_2(C) - \lambda_2(x) \lambda_3(C) + \ldots \quad \textrm{in}
\quad K.\]

{\bf Proof.} We prove Lemma 6.2 by induction on $d$. The case $d=1$ is a special
case of Lemma 4.1(d). So, we may assume that $d\ge 2$ and that Lemma 6.2 holds for
$d-1$. By the splitting principle (see Theorem 6.1 on p.\ 266 in \cite{AT}), we 
may assume that $C=C'+L$ where $C'\in K$ is of $\lambda$-degree $d-1$ and $L \in K$
is of $\lambda$-degree $1$. Using the induction hypothesis, we obtain:
\begin{eqnarray*}
\lefteqn{\gs_2(\lambda_{-1}(C)x) = \gs_2(\lambda_{-1}(C')(1-L) x) =}\\
&& = (\gs_2((1-L)x) - \lambda_2((1-L)x)C' + \gs_2((1-L)x)\lambda_2(C')- \ldots)
\lambda_{-1}(C').
\end{eqnarray*}
Since $\gs_2((1-L)x) = (\gs_2(x) - \lambda_2(x)L)(1-L)$ and $\lambda_2((1-L)x) =
(\lambda_2(x) - \gs_2(x)L)(1-L)$, this implies:
\begin{eqnarray*}
\lefteqn{\gs_2(C,x)}\\
&& = \gs_2(x) - \lambda_2(x)L - \lambda_2(x)C' + \gs_2(x)LC'+
\gs_2(x)\lambda_2(C') - \lambda_2(x)L\lambda_2(C') - \ldots\\
&& = \gs_2(x) - \lambda_2(x)C + \gs_2(x)\lambda_2(C) - \lambda_2(x) \lambda_3(C) +
\ldots
\end{eqnarray*}
as was to be shown.

Now, let $R$ be a commutative ring and $I$ an ideal in $R$ which is locally 
generated by a regular sequence of length $d$. For any f. g.\ projective
$R/I$-module $V$, we choose an $R$-projective resolution $P.(V)$ of $V$. 
The following statement is a consequence of Theorem 5.1.

{\bf Corollary 6.3.} Let $2$ be invertible in $R$. Then we have:
\[H_kN\Sym^2\gG(P.(V)) \cong \left\{\begin{array}{ll}
\Sym^2(V) \otimes \gL^k(I/I^2) & \textrm{for } k \textrm{ even}\\
\gL^2(V) \otimes \gL^k(I/I^2) & \textrm{for } k \textrm{ odd}.
\end{array}\right.\]

{\bf Proof.} By the Eilenberg-Zilber theorem, we have:
\[H_kN\Sym^2\gG(P.(V))\cong H_kN(\gG(P.(V))^{\otimes 2})^{\gS_2} \cong
H_k(P.(V)^{\otimes 2})^{\gS_2}.\]
Therefore, Corollary 6.3 follows from Theorem 5.1.

{\bf Theorem 6.4.} Let $d=2$. Then we have:
\[H_kN\Sym^2\gG(P.(V)) \cong \left\{ \begin{array}{ll}
\Sym^2(V) & \textrm{for } k=0\\
\gL^2(V) \otimes I/I^2 & \textrm{for } k=1\\
D^2(V) \otimes \gL^2(I/I^2) & \textrm{for } k=2\\
0 & \textrm{for } k \ge 3.
\end{array}\right. \]

{\bf Proof.} Let $F_k$ denote the functor
\[F_k: \cP_{R/I} \ra (R\textrm{-modules}), \quad V \mapsto H_kN\Sym^2\gG(P.(V)).\]
Obviously, the functor $F_k$ is of degree $\le 2$. By the Eilenberg-Zilber theorem
(see Application 8.5.3 on p.\ 277 in \cite{We}), we have for all $V,W \in \cP_{R/I}$:
\begin{eqnarray*}
\lefteqn{\rmcr_2(F_k)(V,W) \cong H_kN(\gG(P.(V)) \otimes \gG(P.(W))) } \hspace{2cm}\\
&& \cong H_k \Tot(P.(V) \otimes P.(W)) \cong \Tor_k^R(V,W).
\end{eqnarray*}
Furthermore, we have $\Tor_k^R(R/I,R/I) \cong \gL^k(I/I^2)$ by Example 5.2.
Therefore, by Proposition 1.8, Corollary 1.6, Proposition 1.9, and Example 1.7,
it suffices to show: The plus map $+_2: \rmcr_2(F_0)(R/I, R/I) \ra F_0(R/I)$ is
bijective, $F_1(R/I)$ vanishes, the diagonal map $\gD_2: F_2(R/I) \ra \rmcr_2(R/I,R/I)$
is bijective, and $F_k(R/I)$ vanishes for $k\ge 3$. For this, we may assume that
the ideal $I$ is globally generated by a regular sequence, say $I=(f,g)$. Let $K.$
denote the complex $\ldots \ra 0 \ra R \,\, \stackrel{f}{\ra} \,\, R$ and $L.$ the
complex $\ldots \ra 0 \ra R \,\, \stackrel{g}{\ra}\,\, R$. Then,
$\Tot(K.\otimes L.) \cong \Kos^2(R\oplus R \,\, \stackrel{(f,g)}{\longrightarrow} R)$
is a resolution of $R/I$. For any bisimplicial $R$-module $M..$, let $\gD M..$ denote
the simplicial $R$-module $\gD^\op \ra (R\textrm{-modules})$, $[n] \mapsto M_{nn}$.
The Dold-Kan correspondence (see Proposition 2.1) and the Eilenberg-Zilber theorem
yield natural isomorphisms
\begin{eqnarray*}
\lefteqn{F_k(R/I) \cong H_kN\Sym^2\gG\Tot(K. \otimes L.) \cong
H_kN\Sym^2\gG\Tot(N\gG(K.) \otimes N\gG(L.))}\\
&&\cong H_kN \Sym^2 \gG N \gD(\gG(K.) \otimes \gG(L.)) \cong
H_kN \Sym^2\gD(\gG(K.) \otimes \gG(L.))\\
&& \cong H_kN \gD \Sym^2(\gG(K.) \otimes \gG(L.)).
\end{eqnarray*}
For any two f.\ g.\ projective $R$-modules $P,Q$, the following sequence is exact:
\[0  \ra  \gL^2(P) \otimes \gL^2(Q)  \ra \Sym^2(P\otimes Q)  \ra 
\Sym^2(P) \otimes \Sym^2(Q)  \ra  0 \]
(see p.\ 242 in \cite{ABW}); here, the first map is given by 
\[p_1 \wedge p_2 \otimes q_1 \wedge q_2  \mapsto 
(p_1 \otimes q_1)(p_2 \otimes q_2)- (p_1 \otimes q_2)(p_2 \otimes q_1),\] 
and the second map is given by 
$(p_1\otimes q_1)(p_2 \otimes q_2) \mapsto p_1p_2\otimes q_1q_2$.
Thus, we obtain the following short exact sequence
of bisimplicial $R$-modules:
\[\xymatrix{
\gL^2\gG(K.)\otimes \gL^2\gG(L.) \ar@{^{(}->}[r] & \Sym^2(\gG(K.) \otimes \gG(L.)) 
\ar@{->>}[r] & \Sym^2\gG(K.) \otimes \Sym^2\gG(L.) ;
}\]
applying $N\gD$ to this sequence then yields a short exact sequence of complexes.
By the Eilenberg-Zilber theorem and Proposition 2.4, we have:
\begin{eqnarray*}
\lefteqn{H_kN\gD(\Sym^2\gG(K.) \otimes \Sym^2\gG(L.))}\\
&& \cong H_k\Tot(N\Sym^2\gG(K.) \otimes N\Sym^2\gG(L.))\\
&&\cong H_k\Tot (\Kos^2(f) \otimes \Kos^2(g))\\ 
&& \cong H_k\Tot(K.\otimes L.)\\
&& \cong \left\{\begin{array}{ll}
R/I & \textrm{for } k=0\\
0 & \textrm{for } k\ge 1
\end{array}\right.
\end{eqnarray*}
Similarly, we obtain (cf.\ Remark 2.7):
\[H_kN\gD(\gL^2\gG(K.) \otimes \gL^2\gG(L.)) \cong
\left\{\begin{array}{ll}
R/I & \textrm{for } k=2\\
0 & \textrm{for } k\not= 2.
\end{array}\right.\]
This finally implies:
\[F_k(R/I) \cong \left\{\begin{array}{ll}
R/I& \textrm{for } k=0,2\\
0 & \textrm{for } k=1 \textrm{ or } k\ge 3.
\end{array}\right.\]
In particular, this finishes the proof of Theorem 6.4 in the 
cases $k=1$ and $k\ge 3$. Next, we 
prove that the diagonal map $\gD_2: F_2(R/I) \ra \rmcr_2(F_2)(R/I,R/I)$ is bijective.
As above, we have an isomorphism
\[\rmcr_2(F_2)(R/I,R/I) \cong H_2N\left(\gG(\Tot(K. \otimes L.))^{\otimes 2}\right) 
\cong H_2N\gD\left((\gG(K.) \otimes \gG(L.))^{\otimes 2}\right)\]
such that the following diagram commutes:
\[\xymatrix{
F_2(R/I) \ar[rr]^-{\gD_2} \ar@{=}[d] && \rmcr_2(F_2)(R/I,R/I) \ar@{=}[d]\\
H_2N\gD\Sym^2(\gG(K.) \otimes \gG(L.)) \ar[rr]^-{H_2N\gD(\gD_2)} &&
H_2N\gD\left((\gG(K.) \otimes \gG(L.))^{\otimes 2}\right). 
}\]
Furthermore, we have the commutative diagram of bisimplicial $R$-modules
\[\xymatrix{
\gL^2\gG(K.) \otimes \gL^2\gG(L.) \ar@{^{(}->}[r] \ar[d]^-{\id\otimes \gD_2} &
\Sym^2(\gG(K.) \otimes \gG(L.)) \ar@{->>}[r] \ar[d]^-{\gD_2} &
\Sym^2\gG(K.) \otimes \Sym^2\gG(L.) \ar[d]^-{\id\otimes\gD_2} \\
\gL^2\gG(K.) \otimes \gG(L.)^{\otimes 2} \ar@{^{(}->}[r]&
(\gG(K.) \otimes \gG(L.))^{\otimes 2} \ar@{->>}[r] &
\Sym^2\gG(K.) \otimes \gG(L.)^{\otimes 2}
}\]
where the upper row is the exact sequence already constructed above and the lower
row is induced by the natural exact sequence
\[0\ra \gL^2\gG(K.) \ra \gG(K.)^{\otimes 2} \ra \Sym^2\gG(K.) \ra 0.\]
Since both $\rmcr_2(F_2)(R/I,R/I)$ and $F_2(R/I)$ have rank $1$, it therefore
suffices to show that the map
\[H_2N\gD(\id\otimes \gD_2) : H_2N\gD(\gL^2\gG(K.) \otimes \gL^2\gG(L.)) \ra
H_2N\gD(\gL^2\gG(K.) \otimes \gG(L.)^{\otimes 2})\]
is surjective and that $H_2N\gD(\Sym^2\gG(K.) \otimes \gG(L.)^{\otimes 2}) \cong 0$. 
As above, we have an isomorphism
\begin{eqnarray*}
\lefteqn{H_2N\gD(\gL^2\gG(K.) \otimes \gG(L.)^{\otimes 2}) \cong
H_2\Tot(\tilde{\Kos}^2(K.) \otimes L.^{\otimes 2})}\\
&& \cong H_2\Tot (K.[-1] \otimes L.^{\otimes 2})
\cong H_1((R/(f) \,\,\stackrel{\bar{g}}{\ra}R/(f))^{\otimes 2}) \cong
(\bar{g})/(\bar{g})^2 \cong R/I.
\end{eqnarray*}
(Note that $\bar{g}$ is a non-zero-divisor in $R/(f)$.) 
Similarly, we get 
\[H_2N\gD(\Sym^2\gG(K.) \otimes \gG(L.)^{\otimes 2}) \cong 0.\]
The homology class of
the element $(1\otimes 1)\otimes (1\otimes 1) \in (K_1\otimes K_0) \otimes
(L_1\otimes L_0) \cong \Tot(\tilde{\Kos}^2(f) \otimes \tilde{\Kos}^2(g))_2$ obviously
generates $H_2\Tot(\tilde{\Kos}^2(f) \otimes \tilde{\Kos}^2(g))$. The map
$\id \otimes \gD_2$ maps this element to the element $(1\otimes 1) \otimes
(1\otimes 1, -1\otimes 1)$ of 
\[(K_1\otimes K_0)\otimes (L_1 \otimes L_0 \oplus
L_0\otimes L_1) \cong \tilde{\Kos}^2(K.)_1 \otimes (L.^{\otimes 2})_1 \subseteq
\Tot(\tilde{\Kos}^2(K.) \otimes L.^{\otimes 2})_2.\] 
The homology class of
$(1\otimes 1) \otimes (1\otimes 1, -1\otimes 1)$ obviously generates
$H_2\Tot(\tilde{\Kos}^2(K.)\otimes L.^{\otimes 2})$. This finishes the proof
of Theorem 6.4 in the case $k=2$. Similarly, one can show that the plus map
$+_2:\rmcr_2(F_0)(R/I,R/I) \ra F_0(R/I)$ is bijective. Then, the proof of Theorem 6.4 
is complete. However, this final case can also be proved directly as
follows:
\begin{eqnarray*}
\lefteqn{H_0N\Sym^2\gG(P.(V)) \cong \coker(\xymatrix@1{
\Sym^2(P_1) \oplus P_1\otimes P_0 \ar[r]^-{d_0 - d_1} & \Sym^2(P_0)})}
\hspace{3.5cm}\\
&& \cong \coker(\xymatrix@1{P_1\otimes P_0 \ar[r]^-{{\rm can}}&
\Sym^2(P_0)}) \cong \Sym^2(V)
\end{eqnarray*}
(see also Remark 3.3(b)).

{\bf Remark 6.5.} For any functor $F: \cP_R \ra \cP_R$ (with $F(0) =0$) which
commutes with localization, the support of the homology modules
$H_kNF\gG(P.(V))$, $k\ge 0$, is contained in $V(I) = \Spec(R/I)$ since the 
complex $NF\gG(P.(V))|_{\Spec(R)\backslash V(I)}$ is homotopic to the zero complex
(cf.\ section 2). In contrast to Theorem 3.2, Theorem 5.1, and Theorem 6.4,
these homology modules need not to be projective $R/I$-modules, they even
need not to be annihilated by $I$ in general. For instance, if $d=1$, the 
natural short exact sequence of complexes
\[0\ra ND^2\gG(P.(V)) \ra N(\gG(P.(V))^{\otimes 2}) \ra N\gL^2\gG(P.(V))\ra 0\]
together with Remark 3.6 and Example 5.2 yield the long exact sequence
\begin{eqnarray*}
\lefteqn{ 0 \ra H_1ND^2\gG(P.(V)) \ra V^{\otimes 2}\otimes I/I^2 \ra D^2(V)\otimes I/I^2 }\\
&&\ra H_0ND^2\gG(P.(V)) \ra V^{\otimes 2} \ra \gL^2(V)\ra 0,
\end{eqnarray*} 
and one easily sees that $H_0ND^2\gG(P.(R/I)) \cong R/I^2$ if $\textrm{char}(R) =2$.
If $d=2$, the natural short exact sequence of complexes
\[0\ra N\gL^2\gG(P.(V)) \ra N(\gG(P.(V))^{\otimes 2}) \ra N\Sym^2\gG(P.(V))\ra 0\]
together with Theorem 6.4 and Example 5.2 yield the long exact sequence
\begin{eqnarray*}
\lefteqn{0\ra H_2N\gL^2\gG(P.(V)) \ra V^{\otimes 2}\otimes \gL^2(I/I^2) \ra 
D^2(V) \otimes \gL^2(I/I^2) }\\
&& \ra H_1N\gL^2\gG(P.(V)) \ra V^{\otimes 2} \otimes I/I^2 \ra
\gL^2(V)\otimes I/I^2 \\
&& \ra \gL^2(V) \ra V^{\otimes 2} \ra \Sym^2(V) \ra 0.
\end{eqnarray*}
One can show that the map $V^{\otimes 2} \otimes \gL^2(I/I^2) \ra 
D^2(V) \otimes \gL^2(I/I^2)$ is $+_2 \otimes \id$. Hence,
$H_kN\gL^2\gG(P.(V))$ has non-trivial $2$-torsion for $k=1,2$.
Similarly, one may expect that $H_kN\Sym^2\gG(P.(V))$ can have non-trivial $2$-torsion
if $d\ge 3$ and $k\ge 2$. This is also suggested by the calculation of the stable
derived functors of $\Sym^2$ (e.g., see Theorem 4.7 in \cite{JM}).

{\bf Example 6.6.} Let $K$ be a $\lambda$-ring, $C\in K$ of $\lambda$-degree $2$, and
$x\in K$, arbitrary. An easy, but somewhat tedious calculation shows the following
representation of $\gs_3(C,x)$:
\[\gs_3(C,x) = \gs_3(x) - s_1^3(x)C + (\gs_2(x)x\lambda_2(C) + \lambda_3(x)\gs_2(C))
-s_1^3(x)C\lambda_2(C) + \gs_3(x)\lambda_2(C)^2.\]
In view of Theorem 5.1 and Theorem 6.4, this suggests the following conjecture:
If $d=2$, then we have:
\[H_kN\Sym^3\gG(P.(V)) \cong \left\{\begin{array}{ll}
\Sym^3(V) & \textrm{for } k=0\\
L_1^3(V) \otimes I/I^2 & \textrm{for } k=1\\
\tilde{L}_1^3(V)\otimes I/I^2 \otimes \gL^2(I/I^2) & \textrm{for } k=3\\
D^3(V) \otimes \gL^2(I/I^2)^{\otimes 2} & \textrm{for } k=4\\
0 & \textrm{for } k\ge 5,
\end{array}\right.\]
and there is a short exact sequence
\[0\ra D^2(V) \otimes V\otimes \gL^2(I/I^2) \ra H_2N\Sym^3\gG(P.(V)) \ra
\gL^3(V) \otimes \Sym^2(I/I^2) \ra 0.\]

\bigskip

Mathematisches Institut II der Universit\"at Karlsruhe, 
D-76128 Karlsruhe, Germany; {\em e-mail:} 
Bernhard.Koeck@math.uni-karlsruhe.de.

\end{document}